\documentclass[12pts]{amsart}
\usepackage{amsmath}
\usepackage{amsthm}
\usepackage{amssymb}
\usepackage{amscd}
\usepackage{amsfonts}
\usepackage{amsbsy}
\usepackage{enumerate}
\usepackage{graphicx}
\usepackage[dvips]{psfrag}
\usepackage{subfigure}

\usepackage[utf8]{inputenc}
\usepackage[T1]{fontenc}
\usepackage{newunicodechar}

\usepackage[colorlinks]{hyperref}
\usepackage[usenames,dvipsnames]{color}
\usepackage{pdflscape}
\usepackage[all]{xy}

\usepackage{tikz}
\hypersetup{urlcolor=BlueViolet}
\hypersetup{citecolor=blue}

\usepackage{color}

\DeclareGraphicsExtensions{.bmp,.png,.pdf,.jpg,.eps}

\begin{document}

\title[ linear stability]{Equilibrium points and their linear stability  in the planar equilateral restricted four-body problem:
A review and new results}


\author[]{Jos\'e Alejandro Zepeda Ram\'{\i}rez}
\address{Dept.  de F\'{\i}sica, UAM--Iztapalapa,
	09340 Iztapalapa,  Mexico City,  Mexico.}

\author[]{Martha Alvarez-Ram\'{\i}rez}
\address{Dept. de Matem\'aticas, UAM--Iztapalapa,
09340 Iztapalapa,  Mexico City,  Mexico.}

\email{alex.zepeda.r@gmail.com, mar@xanum.uam.mx}

\maketitle

\begin{abstract}
In this work, we  revisit  the  planar restricted four-body problem to study the dynamics of an infinitesimal mass under
the gravitational force produced by three heavy bodies with unequal masses, forming an equilateral triangle configuration.
We unify known  results about  the  existence and linear stability of the equilibrium points of this problem which have been obtained earlier,
either as   relative equilibria or a central configuration of the planar restricted $(3 + 1)$-body problem. 
It is the first attempt in this direction.
A systematic numerical investigation is performed to obtain the resonance curves  in the mass space. 
We use  these curves to  answer the question about the existing boundary between the domains of linear stability and instability.  
The characterization of   the total number of stable points found inside the  stability domain is  discussed. 
\end{abstract}


\section{Introduction}
The Newtonian planar $n$-body problem reads as the study of the dynamics of $n$ point particles
with masses $m_i$ and positions$q_i$, $i=1,\dots , n$, moving according to Newton's law of motion.

The equations of motion of the $n$-body problem are 
\begin{equation}\label{eq_newton}
m_j\ddot{q}_j = \sum_{i\neq j} \frac{m_im_j(q_i-q_j)}{r_{ij}^3}, \qquad 1 \leq j \leq n.
\end{equation}
where $r_{ij}=|q_i-q_j|$  is the Euclidean distance  between $q_i$ and $q_j$, and we have chosen the units of length in order that the gravitational constant be equal to one. Let $q=(q_1,\dots , q_n)$ be the $2n$-dimensional configuration vector of the primary bodies.

In the planar Newtonian $n$-body problem the simplest motions, called {\em homographic} solutions, are such that the configuration is constant up to rotation and scaling. {\em Relative equilibria} are homographic solutions with the property that the system rotates about its center of mass $c\in{\mathbb R}^2$ as a rigid body and its angular velocity $\omega\neq 0$ is constant. In a rotating coordinate system 
they become equilibrium solutions  of the $n$-body problem, hence the name. Such a solution is possible if and only if
the initial positions $q_i(0)$ satisfy the algebraic equations
\begin{equation}\label{eq_newtoncc}
\lambda (q_j - c) = \sum_{i\neq j} \frac{m_im_j(q_i-q_j)}{r_{ij}^3},  \qquad 1 \leq j \leq n,
\end{equation}
for some a positive constant  $\lambda$.

A configuration $q=(q_1,\dots , q_n)\in {\mathbb R}^{2n}$ of the planar $n$-body problem, 
satisfying (\ref{eq_newtoncc}), is called a {\em central configuration}.   As a consequence, a relative equilibrium is a  central configuration that rigidly rotates about its center of mass.
The reader is addressed to Chapter 2 of \cite{Moeckel2015} for a thorough treatment of this topic.

The number of  central configurations of the planar $n$-body problem for an arbitrary given set of positive masses have long been established only for $n=3$. Up to rotations and translations, there are always exactly five 
classes of central configurations for each choice of positive masses: Two of these are three bodies of arbitrary mass located at the vertices of an equilateral triangles (Lagrange's solutions) and the remaining three are collinear central configurations (Euler's solutions).
However, a complete classification is not known for $n > 3$. Even the finiteness of the central configurations is a very difficult question.
For the general four-body problem, the finiteness  of the relative equilibria was settled by  Hampton and Moeckel  \cite{Hampton2006}. They showed that there are at least 32 and at most 8472 such equivalence classes, including the 12 collinear ones.

A special case of the $n$-body problem is the limiting case in which one of the masses tends to zero.
In the  planar restricted  ($N + 1$)-body problem,  one is asked to consider the motion of a particle of infinitesimal
mass moving in the plane under the influence of the gravitational attraction of $N$ finite particles (called {\em primaries}) that move 
around their common center of mass by retaining an orbit solution of the $N$-body problem. The
infinitesimal mass body is supposed to have no gravitational effect on the other $N$ bodies.

In what follows, we will focus on the planar, circular, restricted four-body problem, 
where the three bodies with positive mass are located at the vertices of an
equilateral triangle (Lagrange's central configuration) rotating on circular orbits about their common center of mass.
 We will refer to this problem as the planar, equilateral, restricted four-body problem, hereafter ERFBP.
This problem originates from the work of  Pedersen \cite{Pedersen1}.

The contribution in this paper is twofold. First, we shall give a state-of-the-art review of count, location  and linear stability  of the equilibria in the ERFBP. We then provide  new contributions found by us concerning to the 
boundary between the domains of linear stability and instability of 
the full set of equilibrium points of the ERFBP in the rotating frame.

Before going into this, it is worth noting that one  difficulty in counting the exact number of  equilibrium points (relative equilibrium or central configuration) is due to the bifurcations. Thus, some others interesting papers dealing with the equilibrium set and its bifurcations  should  be  mentioned.

There has been a substantial amount of analytic and numerical work involving the  number of equilibria of the ERFBP.
  In order to  be self-contained and to make this paper easily understandable to the reader, we
  summarize briefly some known results about  equilibrium points  of the ERFBP.

We start with  the work by   Pedersen \cite{Pedersen1}, who   made a combination of numerical and analytical methods
to compute the  number and positions of   equilibrium points   for the infinitesimal mass of the ($1+3$)-body problem,
when the three large masses form a Lagrangian equilateral triangle.
He found that,  there can be 8, 9, or 10  equilibrium positions, depending on the values of the primary masses.
Moreover,  Pedersen proved that the set of degenerate  equilibrium points is a simple closed curve contained in the
interior of the triangle of positive masses, namely in the simplex $\Sigma$ in the masas space (that will appear later in this work). 
Here,  we will denote this curve by  $\mathfrak{B}$. He proved that, on the bifurcation curve $\mathfrak{B}$, there are 9  equilibrium points.
Pedersen's numerical calculations were later confirmed  in a paper due to  Sim\'o \cite{Simo1978}, where a numerical study was done for  the number of 
relative equilibrium solutions in the four-body problem for arbitrary masses.

In \cite{Arenstorf}, Arenstorf outlines some analytical  proofs of the main results contained in \cite{Pedersen1}.  
In the last part of his paper, he emphasizes that a careful mathematical analysis and rigorous calculations required to
prove these results  are contained  in the Ph.D. thesis of his student  Gannaway \cite{gannaway}.
As it turns out, in Gannaway's dissertation, there are only a few  analytical evidences for particular assertions.
Nevertheless, most  of Pedersen's substantial affirmations about central degenerate configurations, bifurcations   
and counting  are verified once again only by a thorough numerical analysis.
These studies were  eventually completed  by Barros and Leandro \cite{leandrob2011, leandrob2014}.
They  were able to give a mathematically rigorous computer-assisted   proof proving  that,
$\mathfrak{B}$ is a simple, closed, continuous curve, which lies inside the triangle  $\Sigma$ formed by the positive masses.
Besides, Leandro and Barros  also confirmed that  there are either 8, 9, or 10 equilibrium
solutions (depending on the primary masses), and proved that   6 of them are  outside of  the Lagrange equilateral triangle formed by the primary bodies.
Recently, Figueras et al. \cite{jordi2022}   gave a new proof to the one performed by   Barros and Leandro  \cite{leandrob2011, leandrob2014}, which is also based on
computer-assisted methods, but they   applied real analysis techniques instead of (complex) algebraic geometry to counting relative equilibria in the ERFBP.
In contrast to the original proof by Leandro and Barros, their proof does not require any difficult computation.

The  finiteness of the number of equilibria (central configurations) in the  ERFBP is demonstrated by Kulevic et al.
\cite{gareth2009}. They used  tools from  algebraic geometry to state that the number of equilibria in the ERFBP is finite for any choice of masses,
 and is bounded above by 196.  However, they claim that most of these solutions of the equilibrium equations found by them, are physically meaningless. The numerical simulations suggest that the true number varies from 8 and 10, depending on the masses. 
These lower estimates are  just as described on  \cite{Pedersen1}, \cite{Arenstorf},  \cite{gannaway}, \cite{Simo1978},  \cite{BalPap2011a},   \cite{leandro2006}, \cite{leandrob2011},  \cite{leandrob2014}, \cite{zotos2020}, \cite{jordi2022}.

On the other side, the paper by  Baltagiannis and  Papadakis \cite{BalPap2011a},  one of the  recent works we consider in this article,  the authors provided an extended list of possible combinations of  primary bodies masses  and their respective number of points of equilibrium. Other works related to ours are Budzco and Prokopeny \cite{BudzkoProkopenya} and Zepeda Ram\'{\i}rez et al. \cite{alejandro1}, where the authors study  the non-linear stability in the case when the mass parameters of the system lie inside  the domain of linear stability points. In the same vein, we will review a recent work by   Zotos \cite{zotos2020} where  previously known results are retrieved.

It should be noted that, in contrast to the restricted three-body problem, the ERFBP with  total mass  normalized to one,  has two parameters masses, and due to this reason the calculations are much  larger and difficult to carry out. 
With the approach followed in the present paper, we are able to verify that  the number  and stability of the equilibria depend on mass parameters of the primary bodies in a continuous way, as seen in previous papers.

This paper is organized as follows. In Section \ref{sec2}, the problem and equations of motion are presented.
In Section \ref{sec3}   the existence and position of the equilibrium points are investigated, while Section  \ref{sec4} 
is devoted to analyze their linear stability.    We quote some of the classic results    and some   recent. 
 The list does not include all the issues but they are rather significant.
A careful  numerical analysis  allow us to conclude that resonance curves play a key role in determining the stability domain, whose border  is given by the $1$:$1$ resonance curve. Actually,  this occurs within and near the three triangular regions in the mass space where the 
Lagrangian relative equilibria are  stable.
The knowledge of  resonance curves leads us to clarify some points that are a little bit obscure in   \cite{zotos2020}. 
In Section \ref{sec5} we outline some conclusions.

Finally, we stress that all our numerical calculations and graphs of the
obtained results have been performed with  the {\sc Mathematica} software. 

\section{ Description of the  problem}\label{sec2}
The Newtonian planar equilateral restricted  four-body problem describes the motion of an infinitesimal particle
under the gravitational attraction of three primaries $m_1$, $m_2$ and $m_3$ arranged in a central configuration of Lagrange, so that the masses are at the vertices of a rotating equilateral triangle, where the rotation has constant angular velocity $\omega$.

Since the equilateral central configuration is possible for all distributions of masses, 
this paper considers the ERFBP where the primaries  are assumed to have arbitrary masses which are normalized
 so that the total mass of the primaries is taken as the unit of mass, that is, $m_1+m_2+m_3=1$. Hence, 
 the number of mass parameters is reduced from three to two.

\begin{figure}[h]
\centering
\includegraphics[scale=0.3]{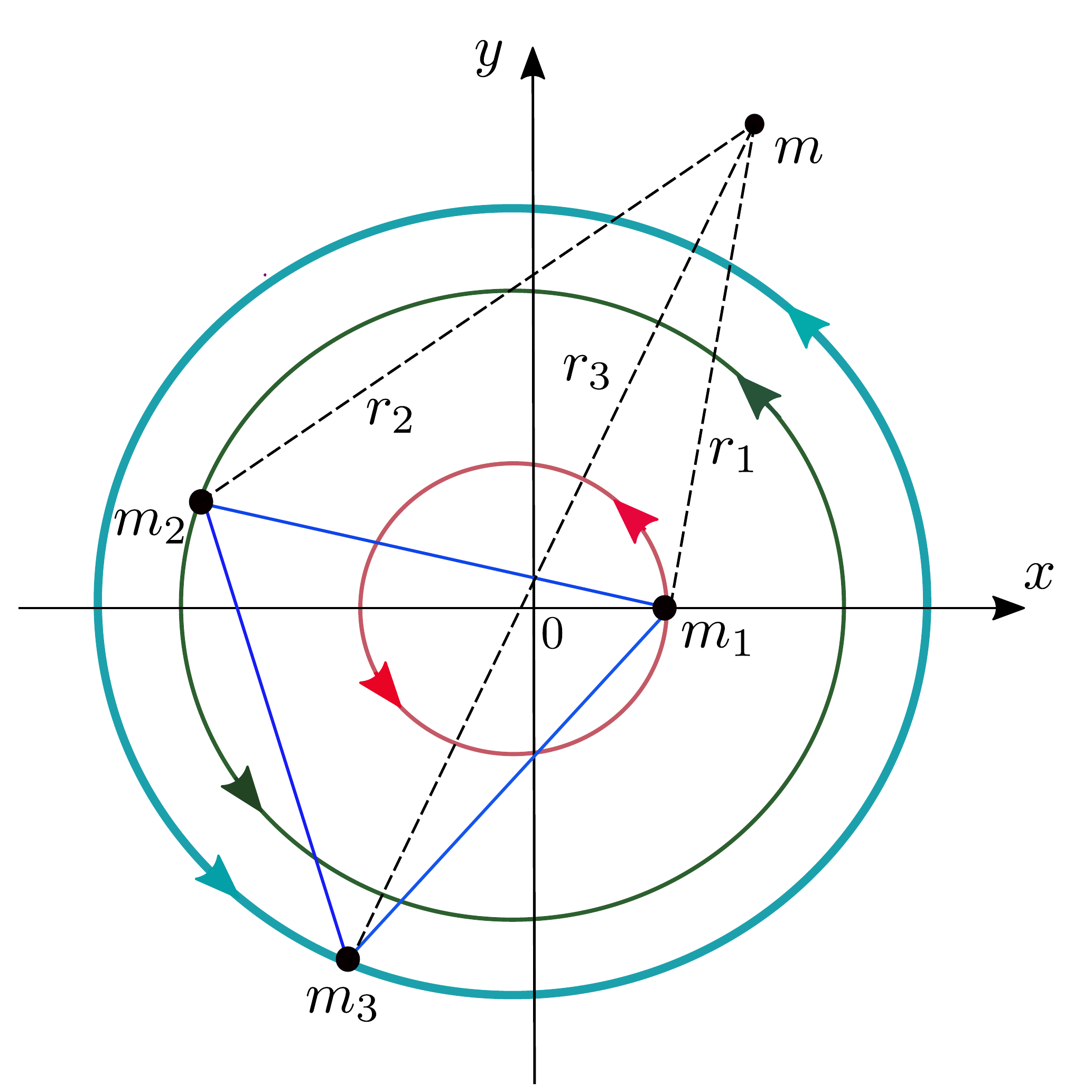}
\caption{The equilateral restricted  four-body problem in  a sidereal frame of reference.}\label{fig1}
\end{figure}

Under the assumption that the center of mass is at the origin of the coordinate system, 
the dynamics of the infinitesimal particle $m$, the system is referred to a rotating (also called synodical) coordinate frame  $(x, y)$,
with uniform angular velocity.
 We orient the triangular configuration so that $ m_1 $ lies on the positive $x$-axis, whose geometry is illustrated in Figure \ref{fig1}.

Define  $K=m_{2}(m_{3}-m_{2})+m_{1}(m_{2}+2m_{3})$. In the synodic reference frame, 
the coordinates  $(x_1, y_1)$,  $(x_2, y_2)$ and $(x_3, y_3)$ of the primaries  are given by
\begin{eqnarray}
	x_{1}=\frac{|K|\sqrt{m^{2}_{2}+m_{2}m_{3}+m^{2}_{3}}}{K}, &\qquad & y_{1}=0, \nonumber\\
	x_{2}=-\frac{|K|\left[(m_{2}-m_{3})m_{3}+m_{1}(2m_{2}+m_{3}) \right] }{2K\sqrt{m^{2}_{2}+m_{2}m_{3}+m^{2}_{3}}},  &\qquad &
	\nonumber y_{2}=\frac{\sqrt{3}}{2}\frac{m_{3}}{m^{3/2}_{2}}\sqrt{\frac{m^{2}_{2}}{m^{2}_{2}+m_{2}m_{3}+m^{2}_{3}}},\\
	x_{3}=-\frac{|K|}{2\sqrt{m^{2}_{2}+m_{2}m_{3}+m^{2}_{3}}}, &\qquad &   y_{3}=-\frac{\sqrt{3}}{2}\frac{1}{m^{1/2}_{2}}\sqrt{\frac{m^{2}_{2}}{m^{2}_{2}+m_{2}m_{3}+m^{2}_{3}}}.\nonumber
\end{eqnarray} 
The equations of motion of the infinitesimal mass $m$   are written as
\begin{equation}\label{four1}
	\begin{array}{l}
\ddot{x}-2\dot{y}=\Omega_{x}, \vspace{0.3cm}\\ 
\ddot{y}+2\dot{x}=\Omega_{y},
	\end{array}
\end{equation} 
where  dots denote derivatives with respect to time $t$, and
\begin{equation}\label{four2}
	\Omega=\Omega(x,y) = \frac{1}{2}(x^2 + y^2) + \frac{m_{1}}{r_1} + \frac{m_{2}}{r_3}+ \frac{m_{3}}{r_2}
\end{equation}
is the potencial function  with $ r_{1}=\sqrt{(x-x_{1})^{2}+y^{2}} $, $ r_{2}=\sqrt{(x-x_{2})^{2}+(y-y_{2})^{2}} $ and  $ r_3=\sqrt{(x-x_{3})^{2}+(y-y_{3})^{2}}$.

The Hamiltonian governing the motion of the infinitesimal particle in these coordinates is
\begin{equation}\label{four5}
	H=\frac{1}{2}(p^{2}_{x}+p^{2}_{y})+yp_{x}-xp_{y} - U(x,y),
\end{equation}
where $p_x=\dot{x}-y$ and $ p_y=\dot{y}+x$  are the conjugate momenta, and 
$$U(x,y)=\frac{m_{1}}{r_{1}} + \frac{m_2}{r_{2}}  +\frac{m_{3}}{r_{2}}$$
is the self-potential.

\begin{figure}[hbt]
\includegraphics[width=0.6\textwidth]{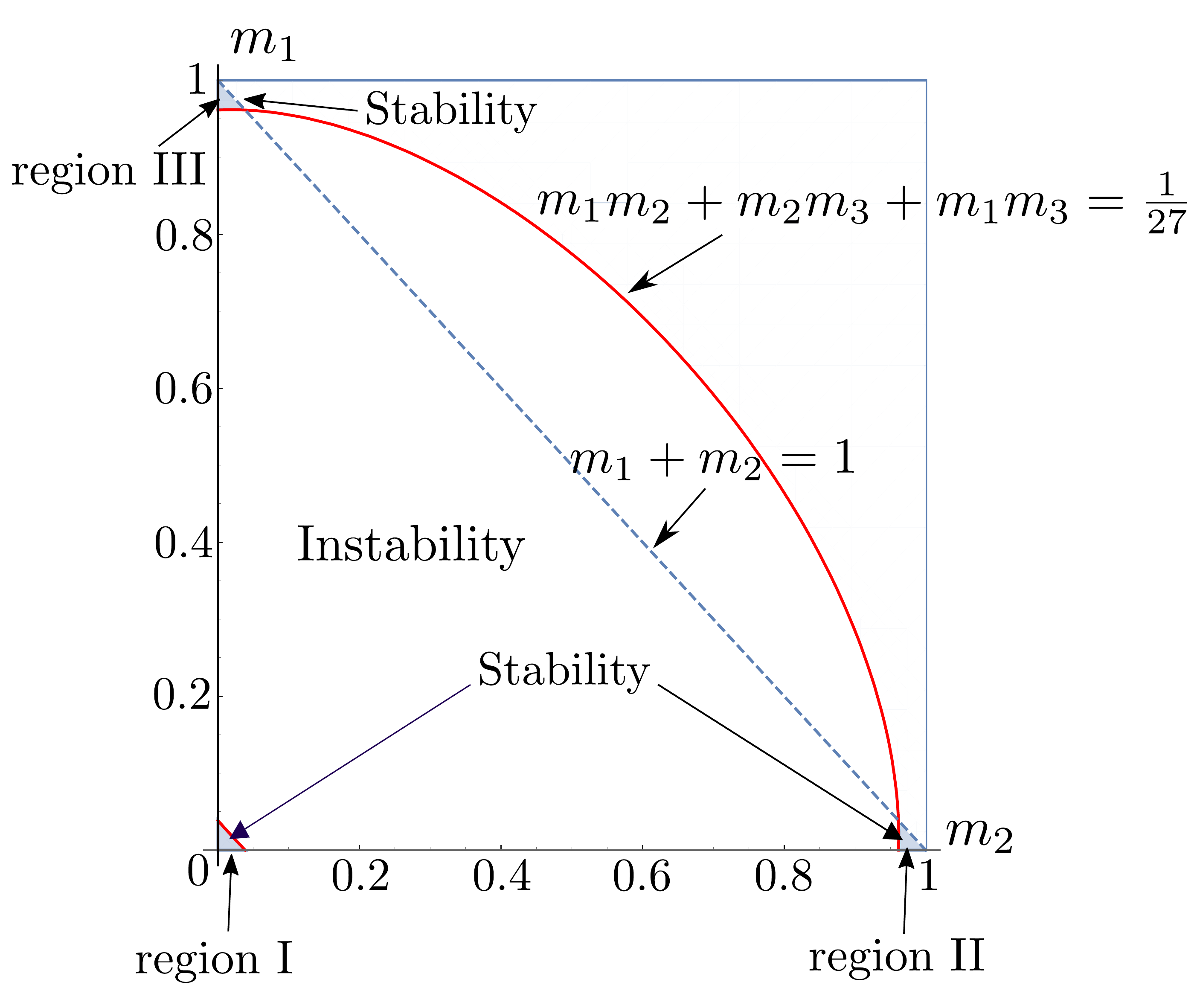}
\caption{The three  small ``triangular" shaded regions I, II and III are  stability domains   of the Lagrange triangle configuration, and white region below the line $m_1+m_2 = 1$ is  instability domain.
The mass parameter of the third primary is $m_3 = 1 - m_1 - m_2$. The red lines correspond to the Routh's critical curve.} \label{fig_reg1}
\end{figure}

At this time the reader should be  warned that  the condition $m_1+m_2+m_3=1$ implies that ERFBP 
depends only on two mass parameters. In particular, we take   $m_3=1-m_1-m_2$ with $ m_1\neq m_2 \neq m_3$, so  the  two 
free parameters will be $m_1$ and $m_2$.  According to Gaschea \cite{Gaschea} and  Routh \cite{Routh},  the triangle configuration
of the primaries is stable   only when  the Routh's stability condition given by inequality
\begin{equation}\label{routh}
 \frac{m_1m_2+m_1m_3+m_2m_3}{(m_1+m_2+m_3)^2} <\frac{1}{27} 
 \end{equation}
is fulfilled. Indeed, this inequality is satisfied only  if one mass of the primaries dominates over  the  other masses.
Figure \ref{fig_reg1}  depicts the stability  regions on the $(m_1,m_2)$ plane, when both inequalities 
$m_1m_2 + m_2m_3 +  m_3m_1< 1/27$ (with $m_3 = 1 - m_1 - m_2$) and $m_1 + m_2 < 1$  are true at the same time;
 the regions of stability are shown in the three gray-shaded areas and instability regions are shown in white, the red curves represent 
solutions to $m_1m_2 + m_2m_3 +  m_3m_1= 1/27$ and  the dashed line represents the condition $m_1+m_2=1$. 
For convenience we label the gray-shaded areas of stability  as follow:  I for the  lower left corner, II for the lower right corner and III for the upper left corner.

We remark that, a simple calculation shows that region II can be   obtained from III by means of reflection with respect to the line $m_1=m_2$.
It follows that  the set of masses on  region II  satisfy  that $m_1$ and $m_3$ are  very small with $m_2$ very large, while on region III 
$m_2$ and $m_3$ are  very small and $m_1$ is very large.

\begin{figure}[hbt]	
\includegraphics[scale=0.3]{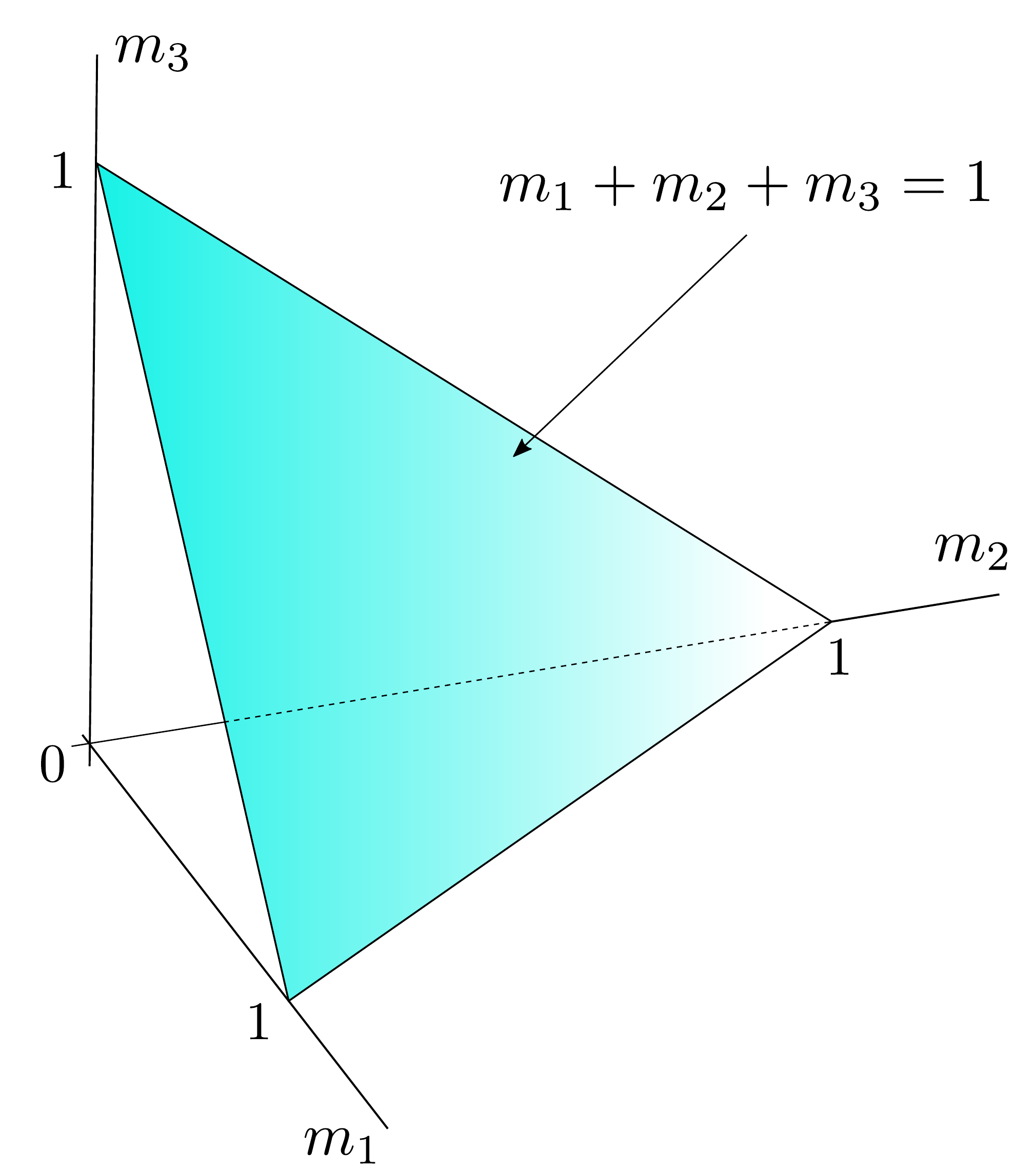}
\caption{Parameter simplex $\Sigma$ in the mass space $(m_1,m_2,m_3)$ normalized so that $m_1+m_2+m_3=1$.  The simplex is formed by the vertices $(1,0,0)$, $(0,1,0)$ and $(0,0,1)$.}\label{masas_map}
\end{figure}

At this stage we need to remind that the normalized mass space can be represented as the 2-simplex 
$$\Sigma = \{(m_1,m_2,m_3)\in \mathbb{R}^3_+  \:  \mid  m_1+m_2+m_3=1, \;  0 \leq m_k \leq 1,\;  k=1,2,3\},$$
 see Figure \ref{masas_map}.  This is an equilateral triangle, whose   edges have length $\sqrt{2}$ and it is called the {\em triangle of masses}.
 It can be seen as  the barycentric coordinates of a point within $\Sigma$ are the masses $m_1$, $m_2$, $m_3$. The sides correspond to mass values of the (2 + 2)-body problem 
(two large and two massless), while the vertices are the masses of the (1 + 3)-body problem.

\section{Equilibrium points}\label{sec3}
We now compute the equilibrium points of  equations (\ref{four1}). By calculating the derivatives of (\ref{four2}) and equaling them to zero the equilibrium points coordinates $(x,y)$ are determined by
\begin{equation} \label{four6}
\begin{array}{l}
\displaystyle{ \frac{\partial \Omega}{\partial x} =  x -\frac{m_{1}(x -x_{1})}{[(x-x_{1})^{2}+y^{2}]^{3/2}}-\frac{m_{2}(x-x_{2})}{[(x-x_{2})^{2}+(y-y_{2})^{2}]^{3/2}}  -\frac{m_{3}(x -x_{3})}{[(x-x_{3})^{2}+(y -y_{3})^{2}]^{3/2}}=0},  \vspace{0.3cm}\\
 \displaystyle{\frac{\partial \Omega}{\partial y} =  y-\frac{m_1 y}{[(x -x_{1})^{2}+y ^{2}]^{3/2}}-\frac{m_{2}(y-y_{2})}{[(x-x_{2})^{2}+(y-y_{2})^{2}]^{3/2}}
-\frac{m_{3}(y-y_{3})}{[(x-x_{3})^{2}+(y -y_{3})^{2}]^{3/2}}=0.}
\end{array}
\end{equation}

\begin{figure}[h!]
	\subfigure[]{	\includegraphics[width=0.42\textwidth]{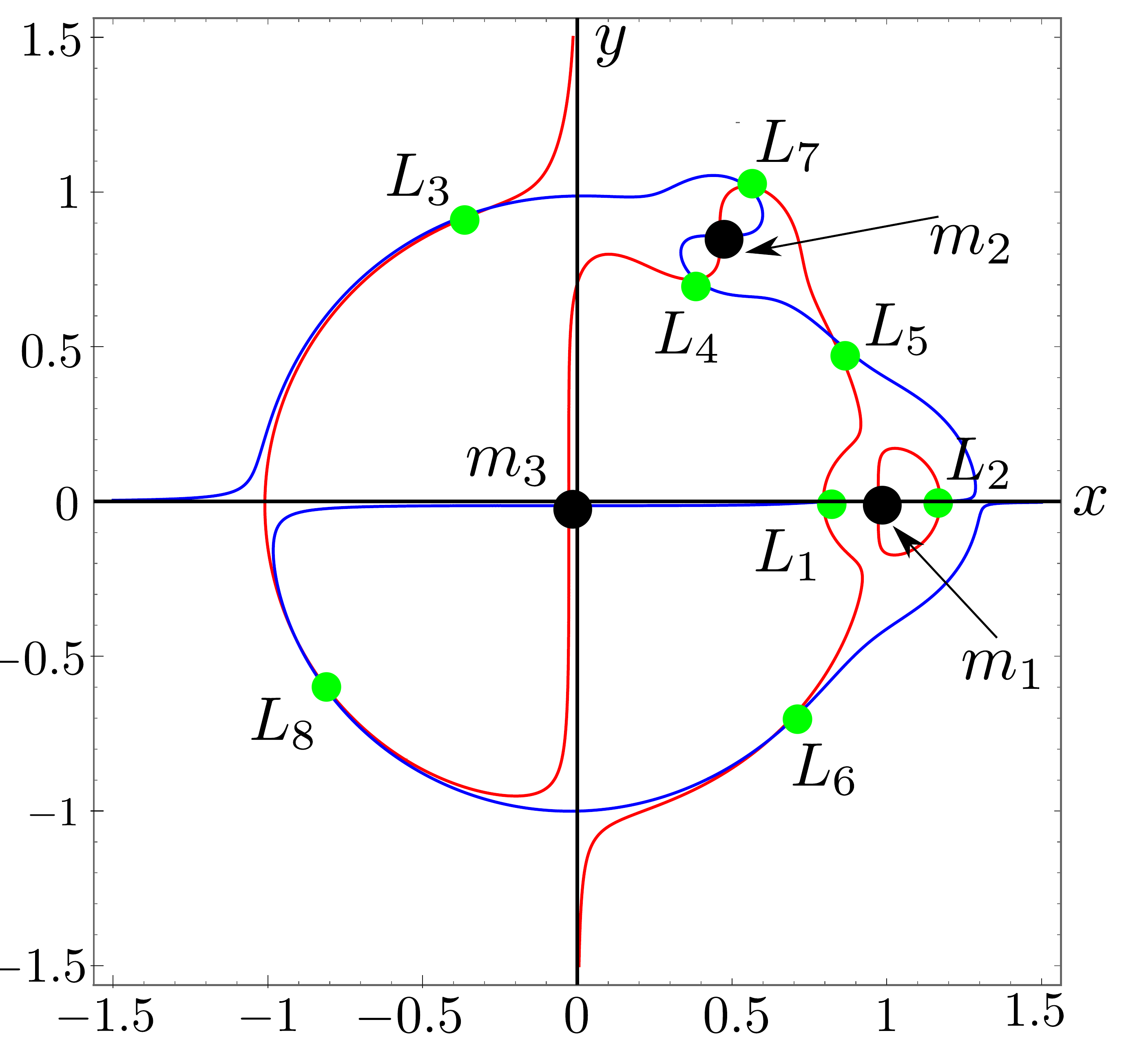}}\label{fig2a}
	\hspace{0.5cm}
	\subfigure[]{	\includegraphics[width=0.42\textwidth]{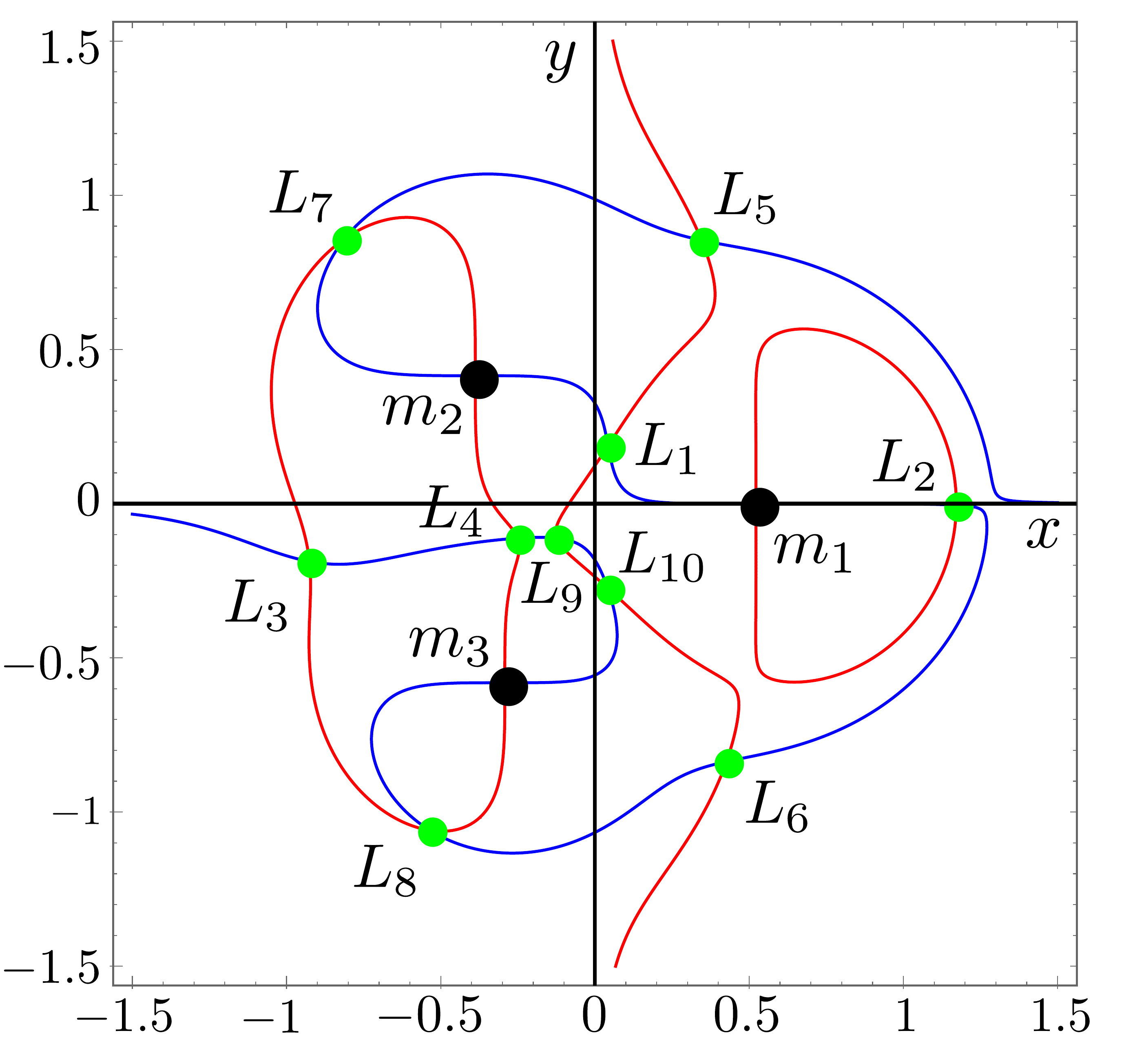}}	
\caption{Equilibrium points in the $x$-$y$ plane located at the intersection 
of the curves $\Omega_x=0$ (red) and $\Omega_y=0$ (blue) for: 
 (a) $ m_1=0.02$ and $ m_2=0.015$   with eight equilibria and (b) $ m_1=0.4$ and $ m_2=0.35$ with ten equilibria. The  green dots denote 
 the position of the equilibrium points and the positions of the primaries are marked by black dots.} \label{fig2}
\end{figure}

Since the system (\ref{four6})  is intractable analytically, the search of equilibrium solutions
 is achieved  by means of numerical methods.
An  intuitive method to locate them relies on finding  intersections of  these zero velocity curves.
We select two sets of mass parameters to plot the positions of the  equilibrium points, through the mutual intersections of the
zero velocity curves shown in Figure \ref{fig2}. 
The labeling for the equilibrium points as $L_i$, $i = 1, 2, \dots, 10$ follows the notation stated in \cite{BalPap2011a}.

\begin{figure}[hbt]	
\includegraphics[scale=0.47]{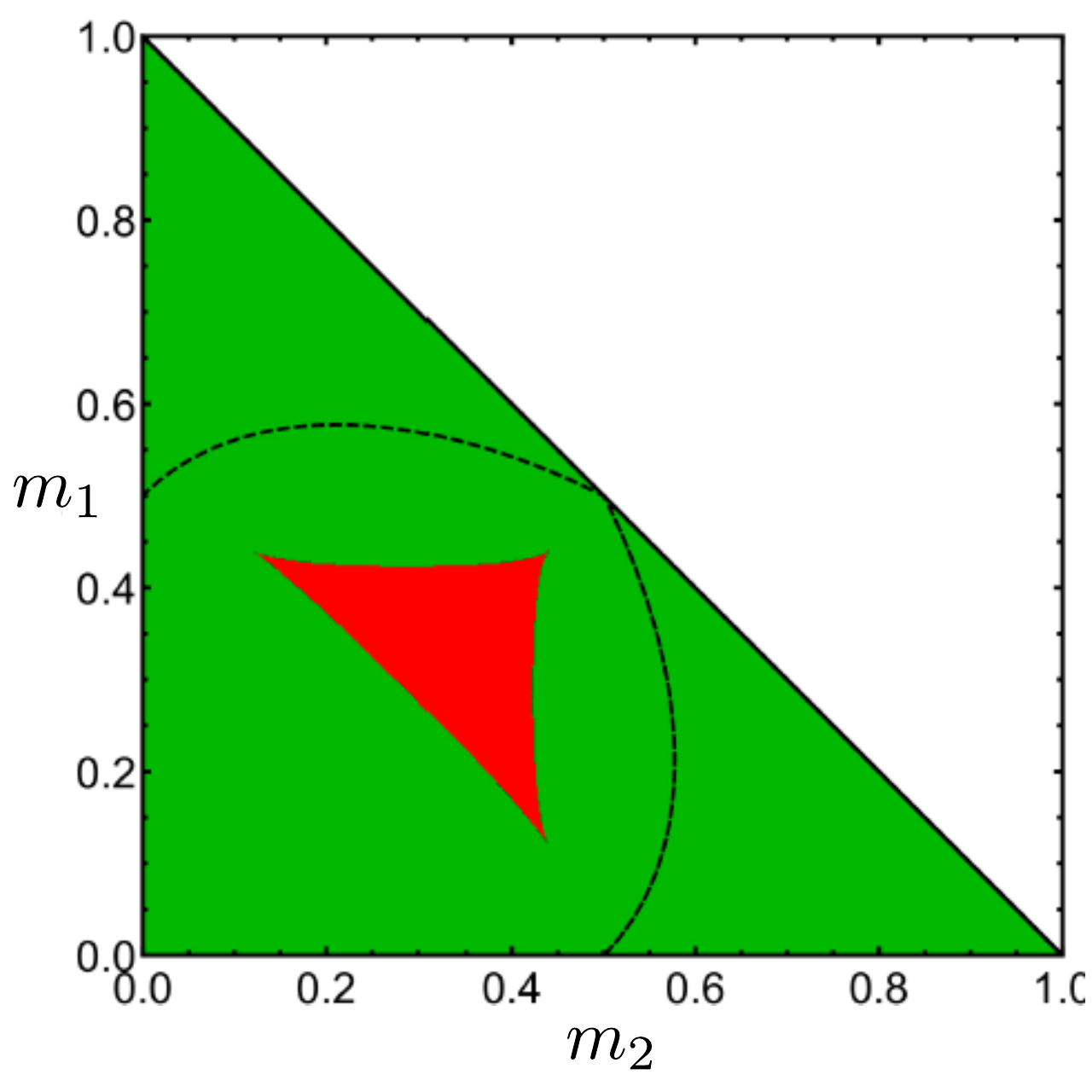}
\caption{Regions on the $(m_1,m_2)$ plane, corresponding to
eight (green) and ten (red) points of equilibrium. 
The black dashed lines are  indicating the set of mass values on which the sign of the relative positions  coordinates of the primaries $m_1$ and $m_2$ is  changing. The border of the red region corresponds to  the bifurcation curve $\mathfrak{B}$.
This  figure is taken from \cite{zotos2020}.}\label{fig_zot1}
\end{figure}

In a recent paper by Zotos \cite{zotos2020},  a detailed numerical study was carried out to show the  number and  location of  equilibrium points of the  ERFBP, for all possible values of the masses of the primaries, as well as its bifurcation set.
Since $m_3 = 1-m_1-m_2$, the results  are  presented   in the simplex defined  by $m_1+m_2 <1$ with $ 0 < m_1, m_2 < 1$  on the  $(m_1,m_2)$ plane,  as in  Figure \ref{fig_zot1}. In fact, Zotos claims that  in the red region there are 10 equilibrium points and in the green region there are only 8. Also, 
the bifurcation curve  $\mathfrak{B}$ is depicted as the border of the central region colored as red, which is almost triangular shaped and cutting the simplex in two components. 
Unfortunately, it seems that Zotos was unaware that the coordinates, number and bifurcation curve $\mathfrak{B}$  of  equilibrium points have been known from earlier works  by, among others,   \cite{Pedersen1},  \cite{Simo1978}, \cite{Arenstorf}, \cite{gannaway},  \cite{leandrob2011}  and \cite{leandrob2014}.

\section{Linear stability of the equilibrium points}\label{sec4}

Once the coordinates of the equilibrium conditions $ (x_0, y_0) $ have been determined, its linear stability can also be studied.
We start by moving the equilibria  to the origin of a coordinate system. Therefore, one sees that the characteristic equation  can be written as:
\begin{equation}\label{four11}
	\lambda^{4}+(4-A_{11}-A_{22})\lambda^2 +A_{11}A_{22}-A^{2}_{12}=0,
\end{equation}
 where 
 \begin{equation}
 \begin{array}{l}
 A_{11} = 1+\displaystyle{\sum_{i=1}^3} \frac{m_i [2(x_0-x_i)^2-(y_0 - y_i)^2]}{[(x_0 - x_i)^2 - (y_0 - y_i)^2]^{5/2}},\vspace{0.3cm}\\
 A_{12} = 3\displaystyle{\sum_{i=1}^3} \frac{m_i [(x_0-x_i)(y_0 - y_i)]}{[(x_0-x_i)^2-(y_0 - y_i)^2]^{5/2}},\vspace{0.3cm}\\
 A_{22} = 1 - \displaystyle{\sum_{i=1}^3} \frac{m_i [(x_0-x_i)^2-2(y_0 - y_i)^2]}{[(x_0 - x_i)^2-(y_0 - y_i)^2]^{5/2}},
  \end{array}  
 \end{equation}
 with $m_1\neq m_2 \neq m_3\neq 0$.
 
By virtue of Lyapunov's theorem  on stability of equilibria for autonomous Hamiltonian systems with two degrees of freedom (see \cite{MeyerOffin}), we have that the equilibria are linearly stable if  (\ref{four11}) has four pure imaginary roots. Indeed, this is  secured by the three following conditions:
\begin{equation}\label{con1s}
\begin{array}{l}
 (4- A_{11}-A_{22})^2 - 4 (A_{11} A_{22}- A_{12}^2) \geq 0, \\
 4- A_{11}-A_{22}>0,\\
 A_{11} A_{22}- A_{12}^2 >0
\end{array}
\end{equation}
which must  be fulfilled simultaneously, whose frequencies  $\omega_1$ and  $\omega_2$  are given by
$$
\omega_{1,2}= \frac{1}{\sqrt{2}}\sqrt{-4 + A_{11}+A_{22}\pm \sqrt{(4- A_{11}-A_{22})^2 - 4 (A_{11} A_{22}- A_{12}^2)}}.
$$

It is known from numerical studies by, among others, Pedersen \cite{Pedersen1}, Arenstorf \cite{ Arenstorf}, Sim\'o \cite{Simo1978} and  Baltagiannis and  Papadakis \cite{BalPap2011a}, 
that the region on the plane $(m_1,m_2)$  where the triangular configuration of the three primaries is stable, there exist eight equilibrium points.
It is noteworthy that  Barros and Leandro \cite{leandrob2014}  used  analytical and computational techniques to prove  that, for all triples $(m_1,m_2,m_3)\in \Sigma$ which are close enough to $\partial\Sigma$, the number of  central configurations is eight. This means we will have eight equilibrium points
 on regions I, II and III.

We note the system of equations (\ref{four6}) is nonlinear, and we  do not know, in advance, where the equilibrium points are located 
for different values of $m_1$  and $m_2$ with $m_3=1-m_1-m_2$.  
For this reason, in order to find the  solutions of (\ref{four6}) as a function of the masses,  we proceed to solve them numerically  in the 
already mentioned region contained in the plane $m_1m_2$. To do so,  we define a dense uniform grid  with a small step, 
which correspond to the initial approximations   in order to apply the command {\em FindRoot}  built within the {\em Mathematica} software.  Then, these solutions, and also  the respective values of the masses,  are inserted into the characteristic equation (\ref{four11}) and thus we derive its linear stability. 

Numerical analysis of conditions  (\ref{con1s}) at each equilibrium point shows that 
 $L_1$, $L_2$, $L_4$, $L_7$ and $L_8 $  are unstable, while  $L_3$, $L_5$ and $L_6$  are the only stable equilibria 
 for any values of the masses from the regions where the Lagrange's configuration is stable. This is in agreement with the results obtained by    Sim\'o \cite{Simo1978},  Budzko and Prokopenya \cite{BudzkoProkopenya},  
 Baltagiannis and  Papadakis \cite{BalPap2011a}  and Zotos \cite{zotos2020}, but some of them used different  labels for the points, see Table \ref{tabb}.
\begin{table}[h] 
\begin{tabular}{|c |c|c|c|}
\hline
Baltagiannis \& Papadakis  & Sim\'o & Budzko \& Prokopenya & Zotos\\
\hline
$L_3$ & $L_4$ & $S_1$ & $L_2$ \\
\hline
$L_5$ & $L_5$ & $S_4$ & $L_3$ \\
\hline
$L_6$ & $L_6$ & $S_7$ & $L_4$ \\
\hline
\end{tabular}
\caption{Labels used to name the linear stable equilibrium points for several authors. We have opted the notation in the same way as  Baltagiannis and  Papadakis \cite{BalPap2011a}.}\label{tabb}
\end{table}

We stress that numerical evidence shows that, for arbitrary masses, 
the equilibrium points $L_3$ and $L_5$ are   located in the upper half plane  ($y>0$), also note that the position of $L_6$  is  in the lower half plane ($y<0$);
see, for instance Figure \ref{fig2a}.

\begin{figure}[hbt]
	\subfigure[Region I,  $B_{L_3}=(0.02012014,0.018114)$.]{	\includegraphics[width=0.42\textwidth]{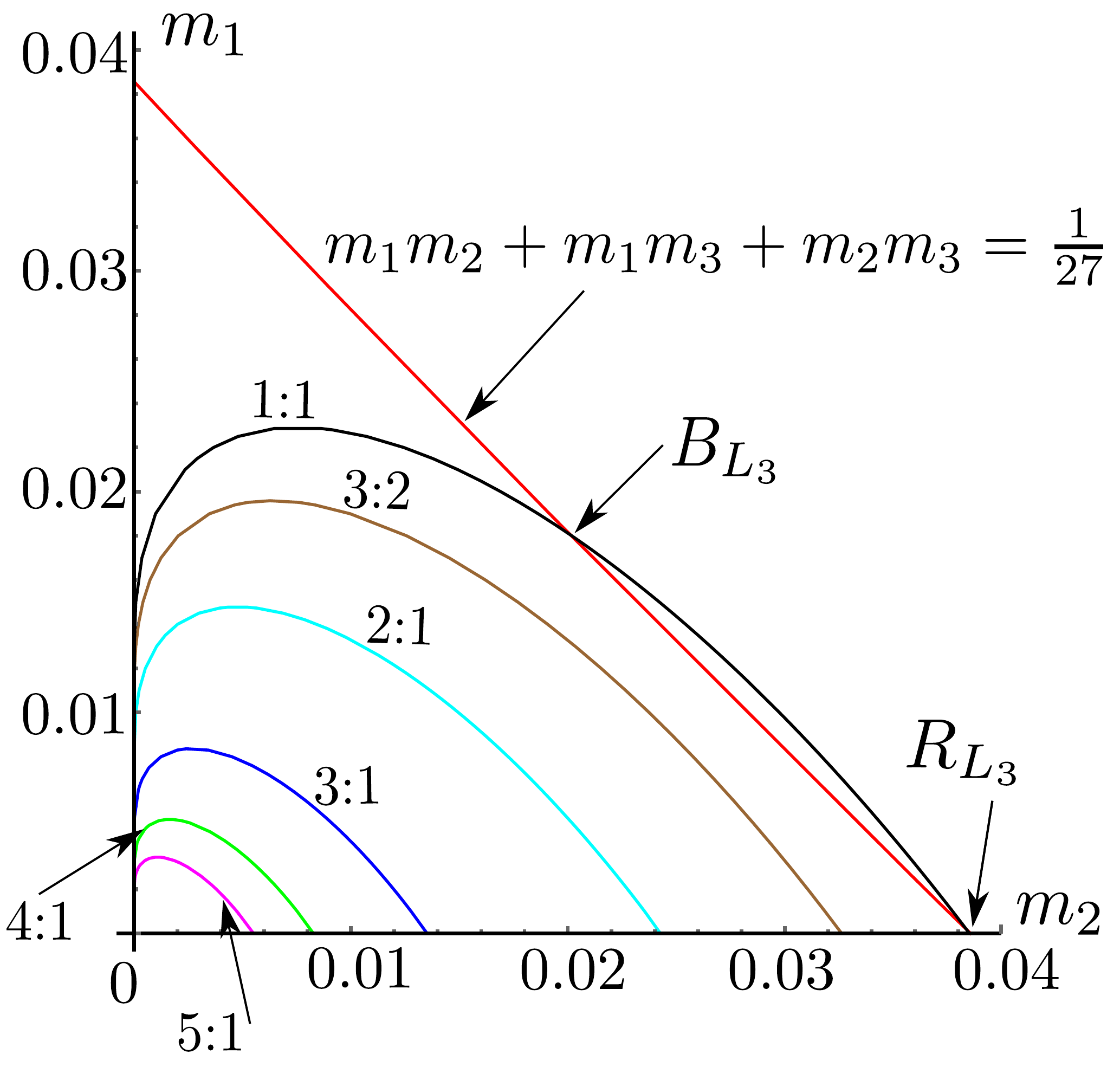}}
	\hspace{0.55cm}
	\subfigure[Region III. ]{\includegraphics[width=0.45\textwidth]{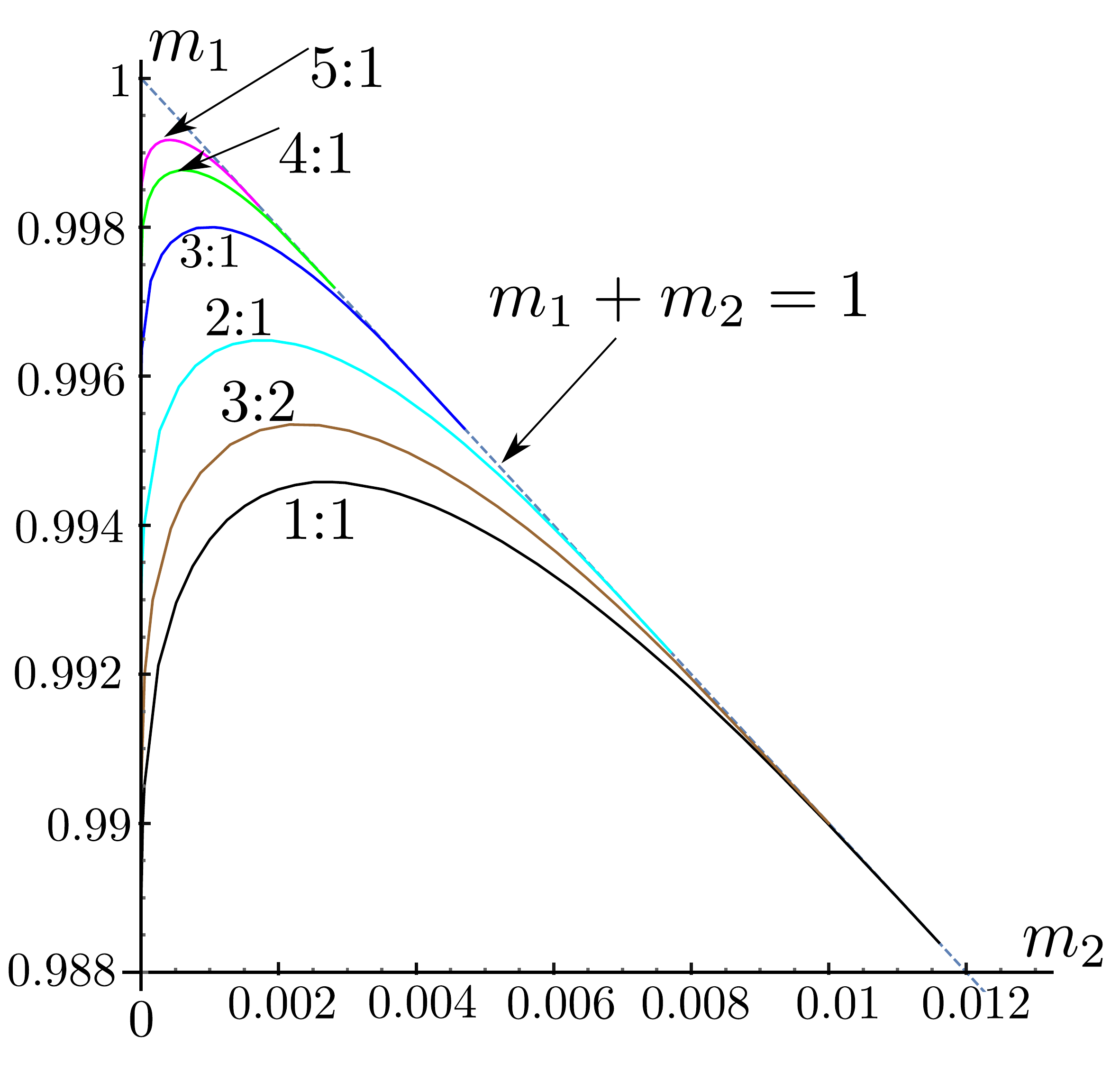}}
	\caption{Resonance curves for the point $L_3$ in regions I and III.}\label{fig4L3}
\end{figure}

\begin{figure}[hbt]
	\subfigure[Region I.]{	\includegraphics[width=0.48\textwidth]{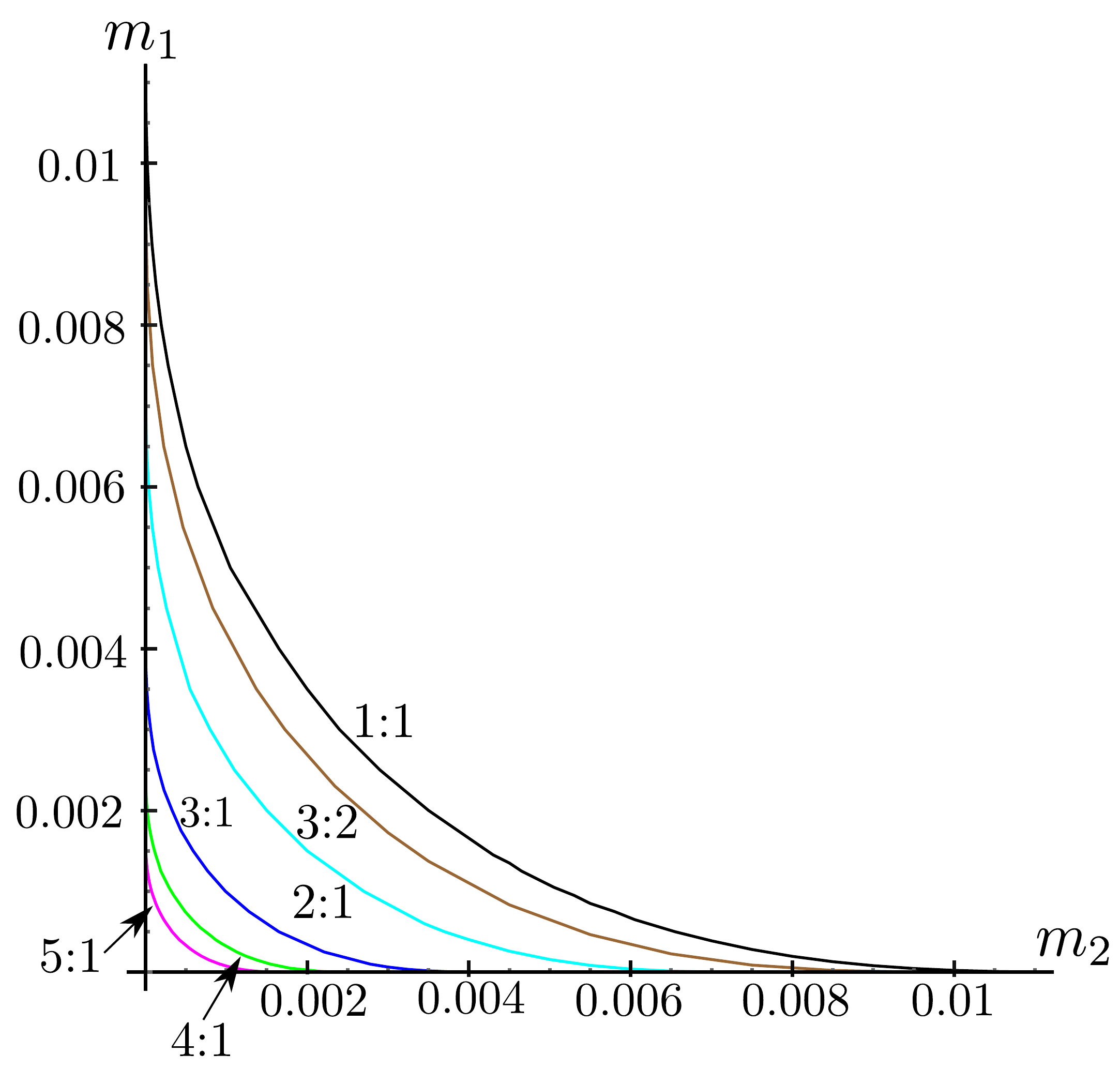}}
	\hspace{0.35cm}
	\subfigure[Region III. $A_{L_5}=(0.02016,0.9619)$  and $R_{L_5}=(0.0385208,0.9614)=(\mu_R,1-\mu_R)$.]{\includegraphics[width=0.45\textwidth]{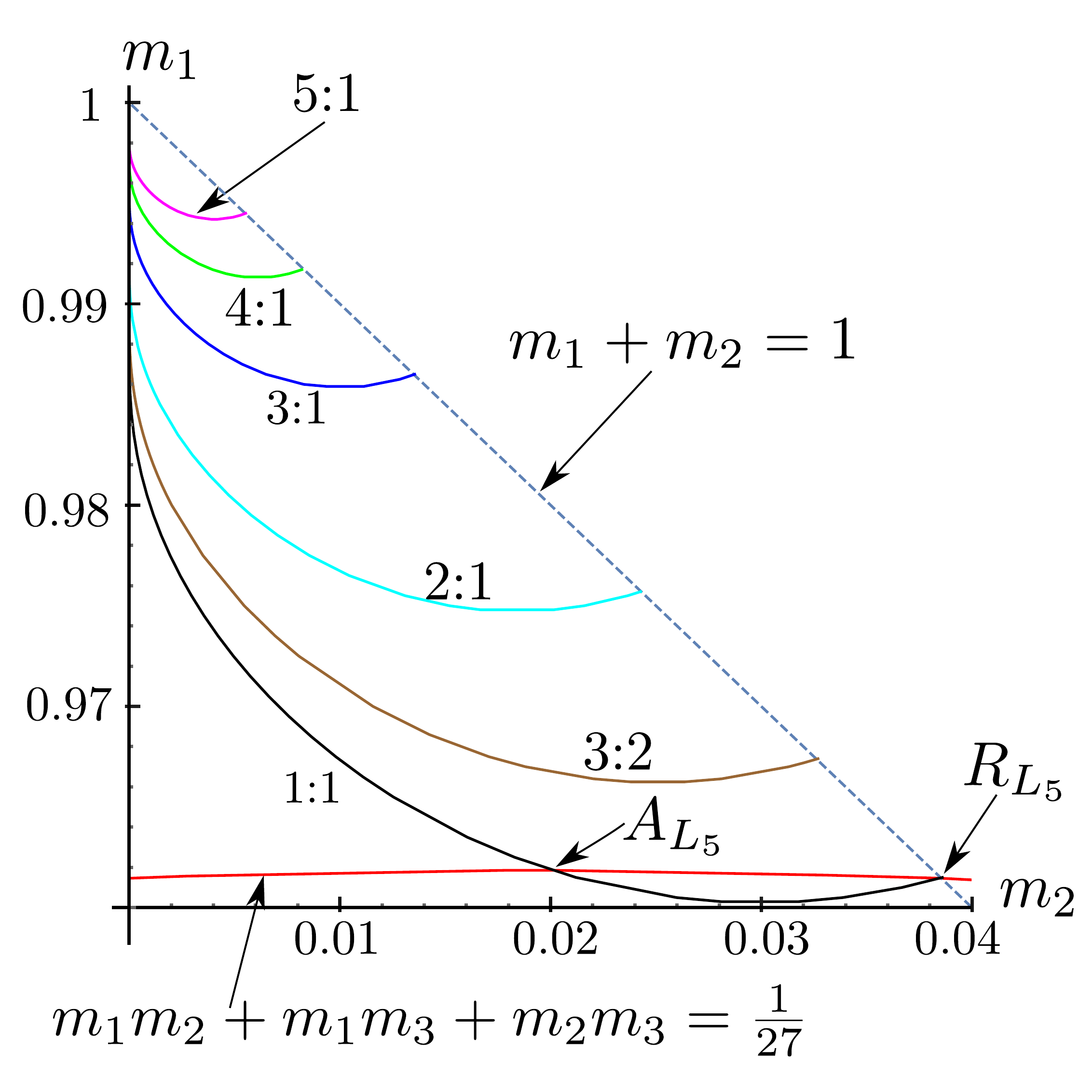}}
	\caption{Resonance curves for equilibrium point  $L_5$ in the regions I and III.} \label{fig4L5}
\end{figure}

\begin{figure}[hbt]
	\subfigure[Region I. $B_{L_6}=(0.018114,0.02012014)$ and $R_{L_6} = (0,0.0385208)=(0,\mu_R)$.]{	\includegraphics[width=0.47\textwidth]{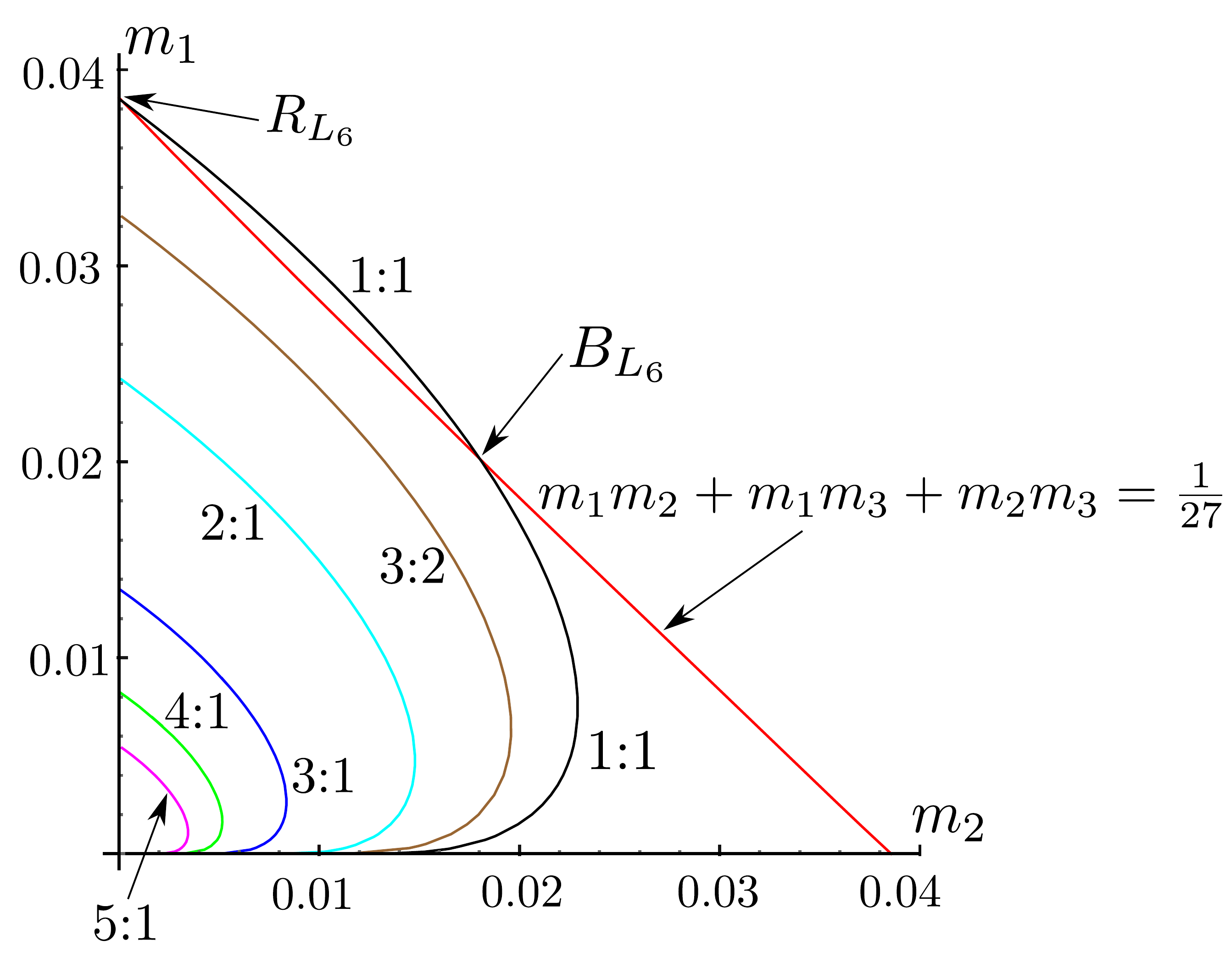}}
	\hspace{0.5cm}
	\subfigure[Region III.  $A_{L_6}=(0.01802,0.9619)$ and ${\mathcal R}_{L_6} =(0.9614,0.0385208)=(1-\mu_R,\mu_R)$.]{\includegraphics[width=0.42\textwidth]{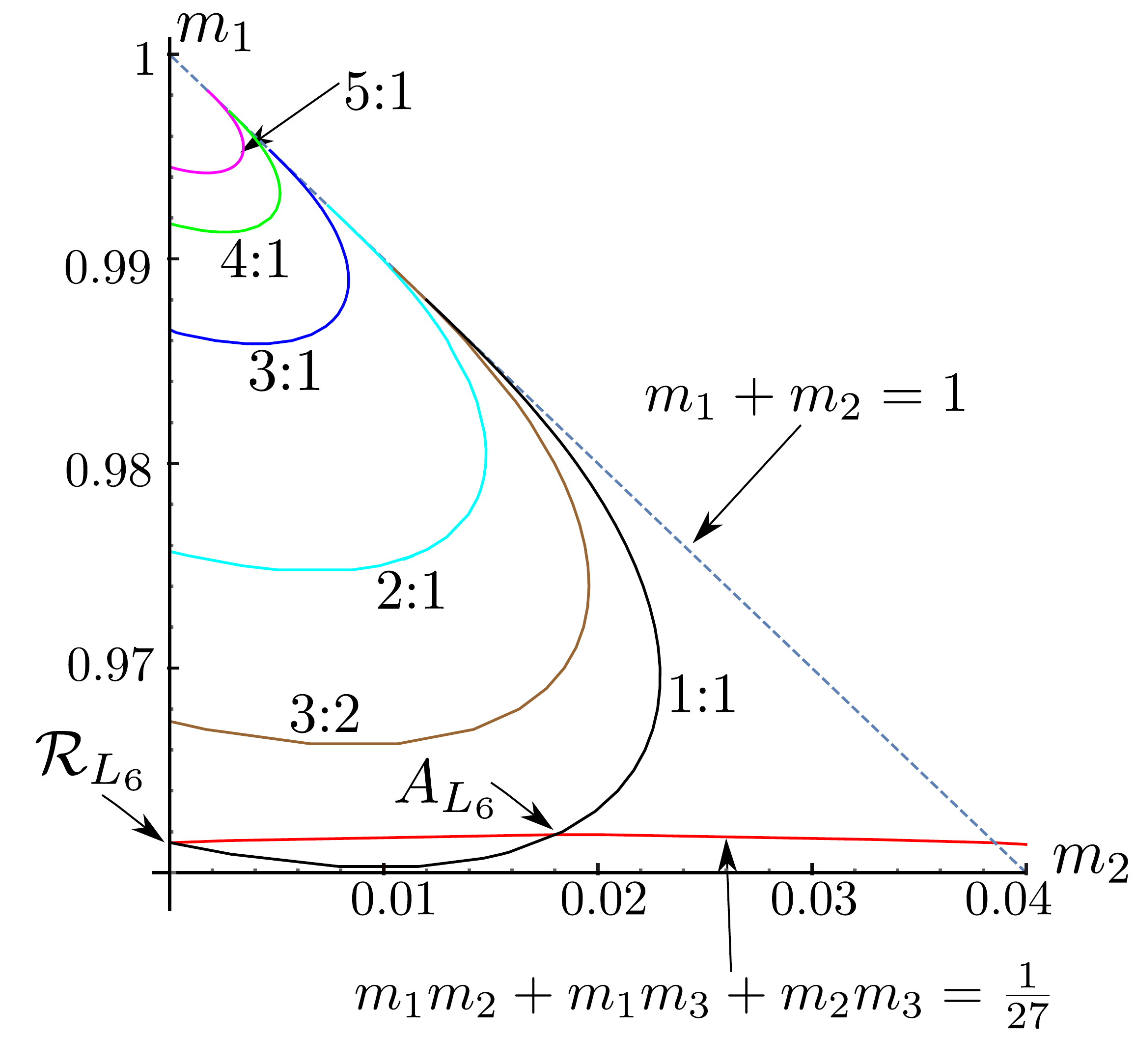}}
	\caption{Graphics of the resonance curves  for equilibrium point  $L_6$ on the regions I and III.} \label{fig4L6}
\end{figure}

Next step is to build numerically  the resonances between frequencies associated to the eigenvalues of the stable
equilibrium points. We state the type and order of resonance  for  each one of the stable equilibria without regard to the stability of the triangle formed by the primaries.  The corresponding   curves are  depicted in  Figures \ref{fig4L3}, \ref{fig4L5} and \ref{fig4L6}, respectively. 

At this stage, we remind  that region II  can be  obtained from region III by means of reflection with respect to the line $m_1=m_2$.
Our computations show that,  the stable equilibrium points  in both regions are $L_3$, $L_5$ and $L_6$. 
In addition, the respective stability domains are symmetric with respect to the line $m_1=m_2$, respectively. Indeed, this is in agreement with the numerical results obtained by  Zotos, which can be checked in the Figures  10(b) and 10(c) in   \cite{zotos2020}.

Now,  we would like to remark that lines $m_2=0$ ($m_3$-axis), $m_3=0$ ($m_2$-axis) and the line $m_1+m_2=1$ on the $m_2m_1$ plane correspond to two copies of the circular restricted three-body problem, respectively. Then, there are points of intersection of
the  resonance curves $1$:$1$, the curve $m_1m_2+m_1m_3+m_2m_3=1$ and   either the axes or  the line $m_1+m_2=1$.
Those points are  $R_{L_3}=(\mu_R,0)$ in Figure \ref{fig4L3} (a), $R_{L_5}=(\mu_R,1-\mu_R)$ in Figure \ref{fig4L5} (b), $R_{L_6}= (0,\mu_R)$ in Figure \ref{fig4L6} (a),
${\mathcal R}_{L_6}= (\mu_R,0)$ in Figure \ref{fig4L6} (b), 
where  $\mu_R = \frac{1}{2} (1-\sqrt{69}/9)\approx 0.038520896504551$ is the Routh's critical mass ratio.

\begin{figure}[hbt]
	\subfigure[Resonance curves in region I.]{	\includegraphics[width=0.44\textwidth]{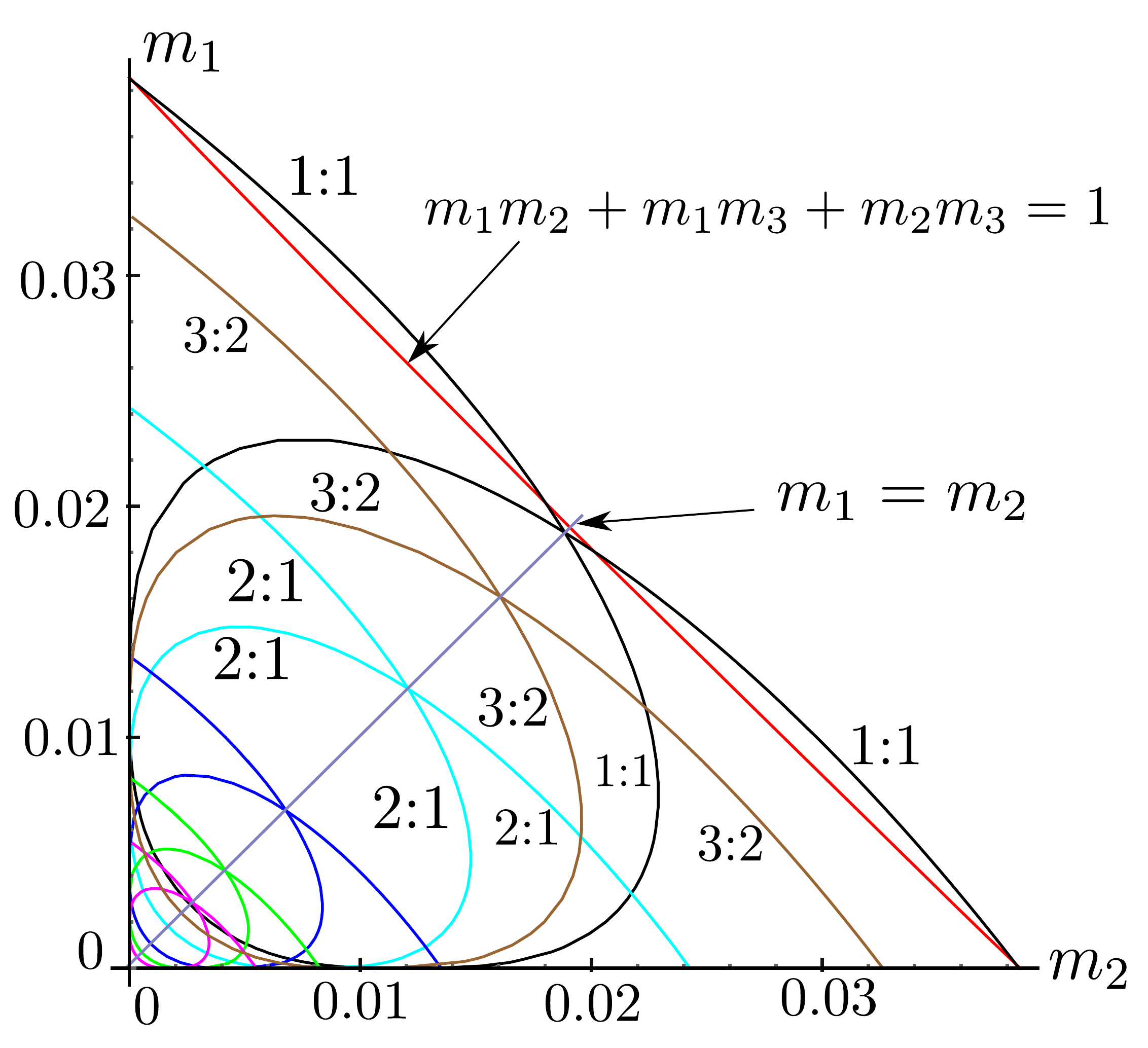}}
	\hspace{0.5cm}
	\subfigure[The enlarged graph.]{\includegraphics[width=0.44\textwidth]{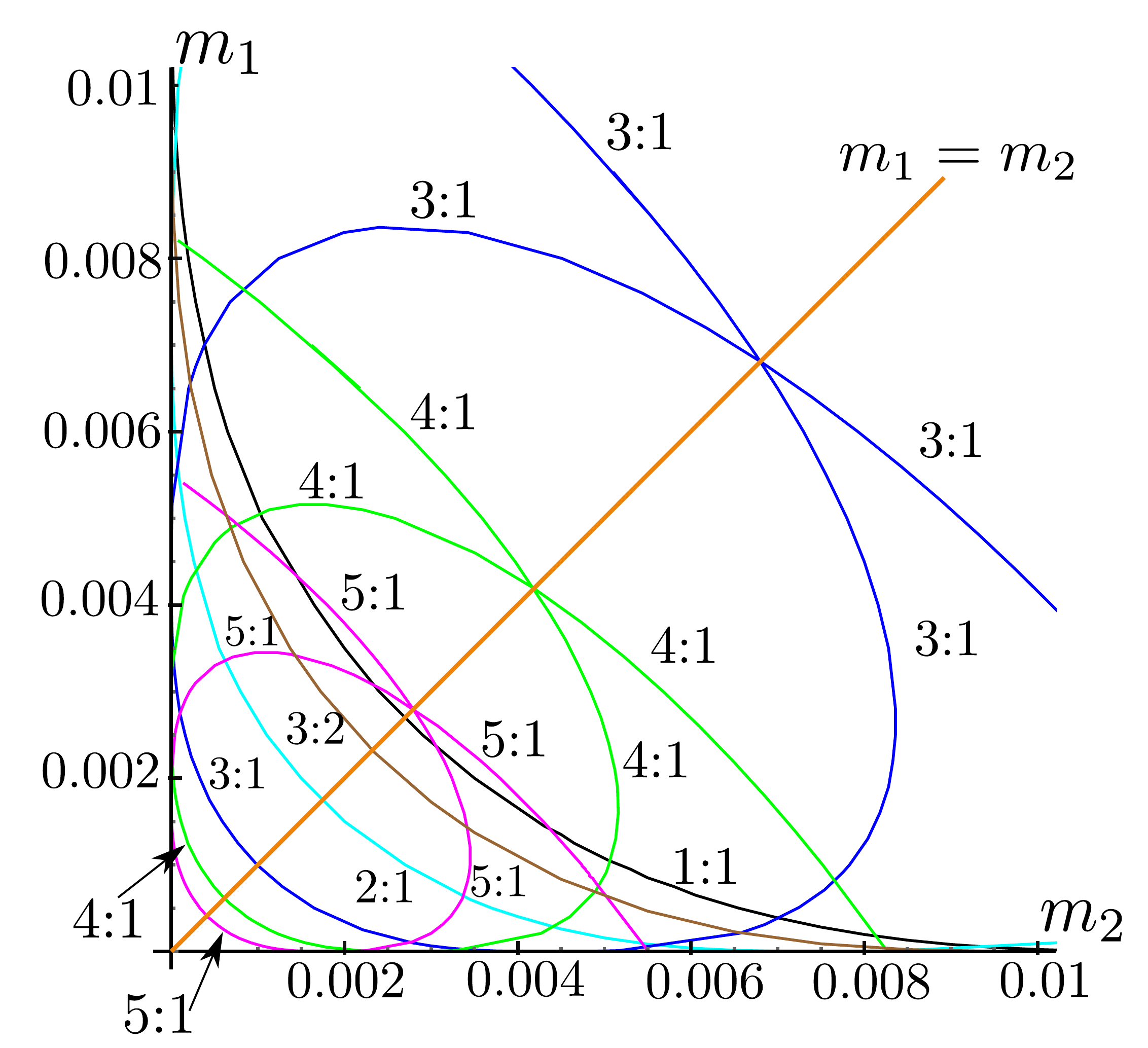}}
	\caption{Plots of all resonance curves   in region I  for equilibrium points  $L_3$, $L_5$ and $L_6$.} \label{fig4LT}
\end{figure}

A striking observation is the fact that the regions where the points of equilibrium are linearly stable 
are bounded  by the  resonance curve $1$:$1$, that is,  the stability domain is given by 
the values of $m_1$ and $m_2$   smaller than their values  on the condition $\omega_1= \omega_2$. 
In a related approach to ours,   Budzco \cite{Budzko2009}
found equilibrium solutions and investigate their linear stability in the ERFBP and  established  that  the stability boundaries are determined by the condition $\omega_1= \omega_2$, but  he only considered the region  I in his study.
Later on, Budzco and  Prokopenya \cite{BudzkoProkopenya}  constructed curves on which the conditions of the third
and fourth order resonances are fulfilled in the region I.

One must note  that,  equilibria $L_3$ and $L_6$ in region I, and also $L_5$ in regions II and III turn out to be linearly stable for some values of parameters $m_1$ and $m_2$   for which  the Routh's condition is not satisfied, then the Lagrange central configuration of the primaries  is unstable.   This fact was observed by  Pedersen  \cite{Pedersen2} who studies the stability of the $(3+1)$-body problem without regard to the stability of the primaries configuration. Also, this very same fact was obtained by  Zotos \cite{zotos2020}.
The bifurcation values found by us are the following:
 for $L_3$ is $B_{L_3}=(0.02012014,0.018114)$,   
for $L_5$ is  $A_{L_5}=(0.02016,0.9619)$, for $L_6$ are $B_{L_6}=(0.018114,0.02012014)$ and  $A_{L_6}=(0.01802,0.9619)$. 
These are given in Figures \ref{fig4L3}(a), \ref{fig4L5}(b),  \ref{fig4L6}(a) and  \ref{fig4L6}(b), respectively.

\begin{figure}[hbt]
	\subfigure[Region III.]{\includegraphics[width=0.54\textwidth]{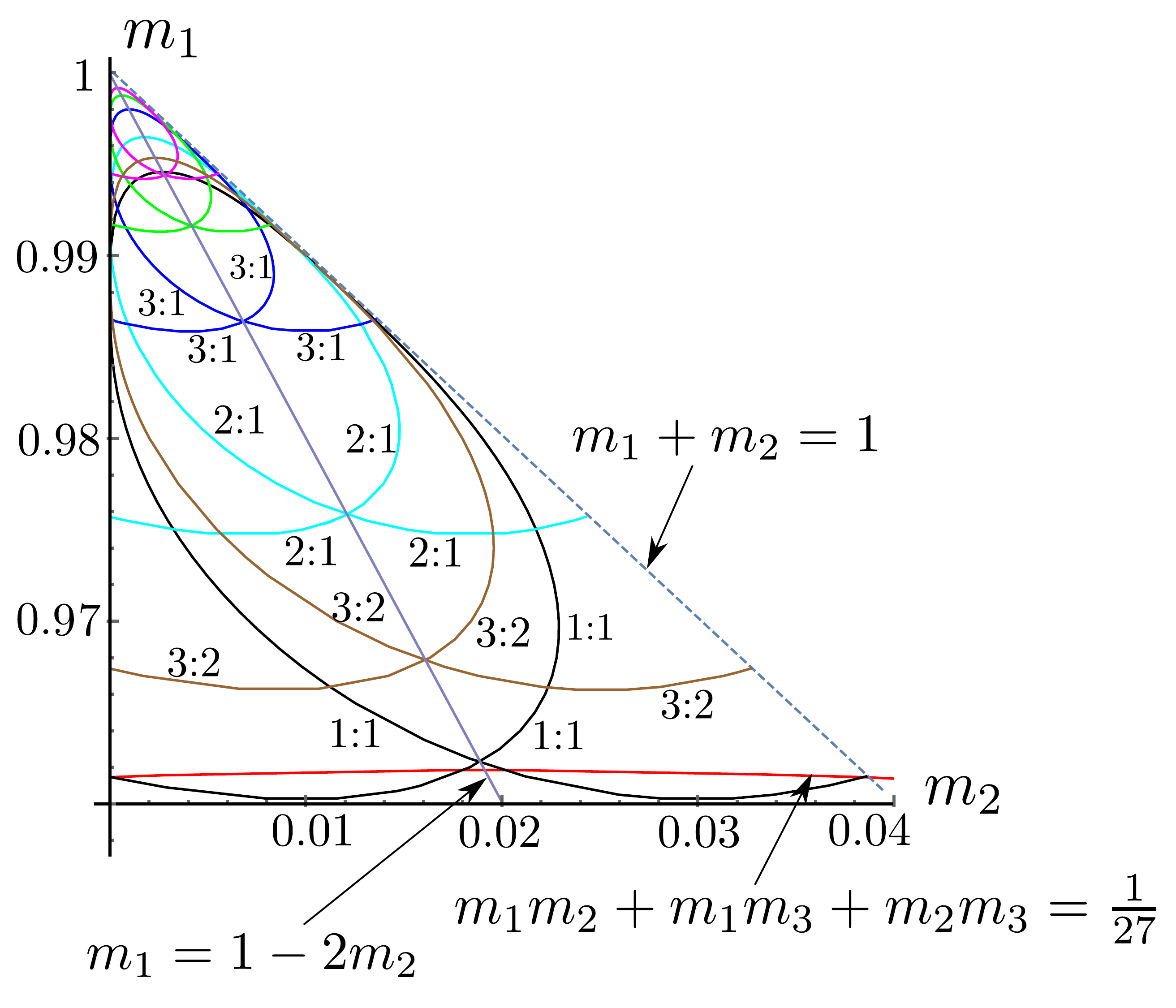}}
	\hspace{0.5cm}
	\subfigure[The enlarged graph.]{\includegraphics[width=0.41\textwidth]{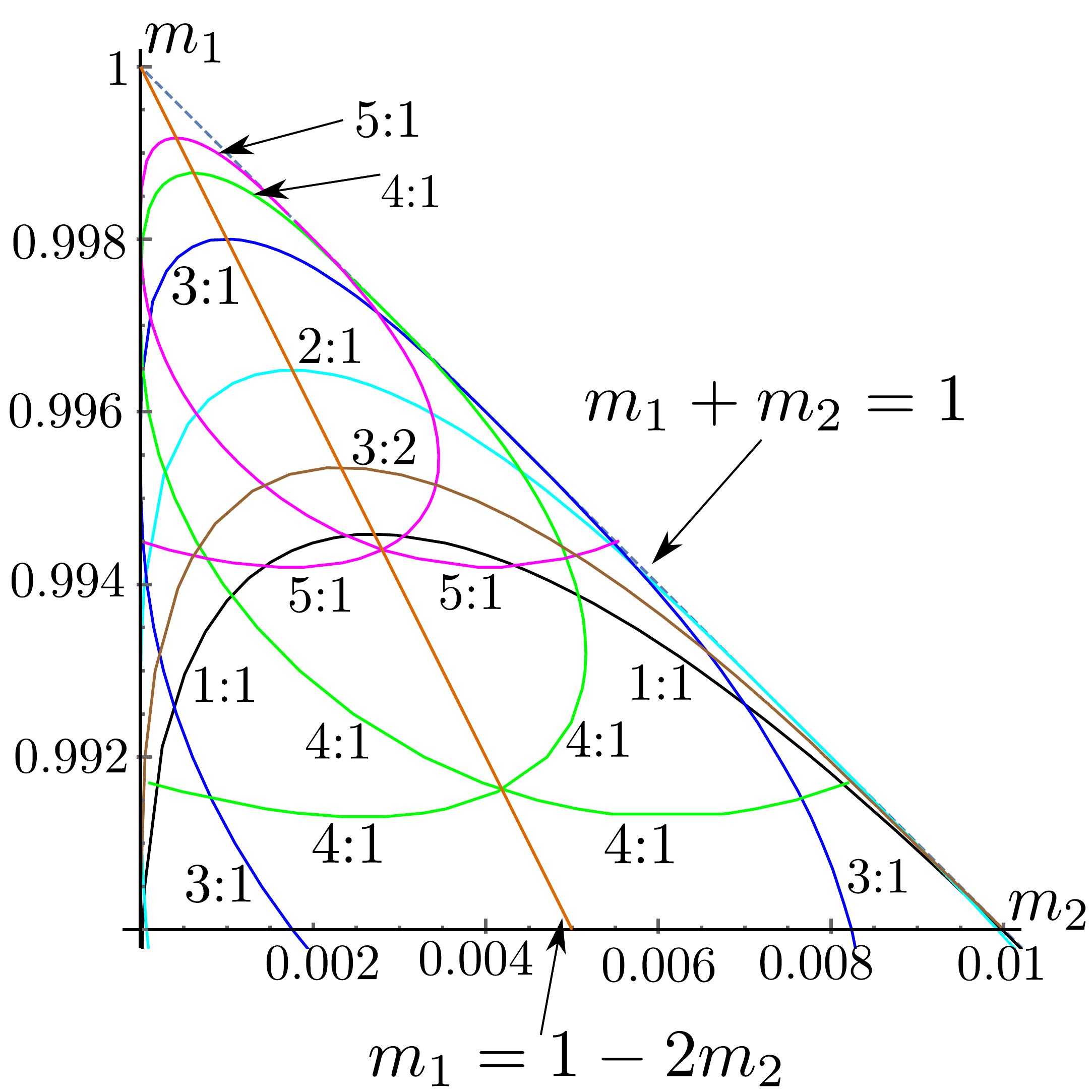}}
	\caption{Resonance curves for equilibrium points  $L_3$, $L_5$ and $L_6$ on the region III. The line  is $m_1=1-2m_2$, that is, $m_2=m_3$.} \label{fig4LTb}
\end{figure}

To determine how many and which stable equilibrium points exist inside or outside  regions I and III,
we draw together the resonance curves of $L_3$, $L_5$ and $L_6$  in a single graph in Figures  \ref{fig4LT} and \ref{fig4LTb}, respectively.

It is not a difficult task to establish a comparison between our  numerical studies and those carried out by  Zotos \cite{zotos2020}.
Figures \ref{fig_cebolla1}(a) and \ref{fig_cebolla2}(a) show the stability domains in  regions I and III, respectively, their topology is exactly that one given by Zotos, but he did not identify the border of the domains, only stated the number of equilibrium points in each subregion and did not say to which point it corresponded. We  go beyond this and are able  to identify how many and what exactly are the respective stable points in the different subregions.

\begin{figure}[hbt]
\centering
\includegraphics[width=0.44\textwidth]{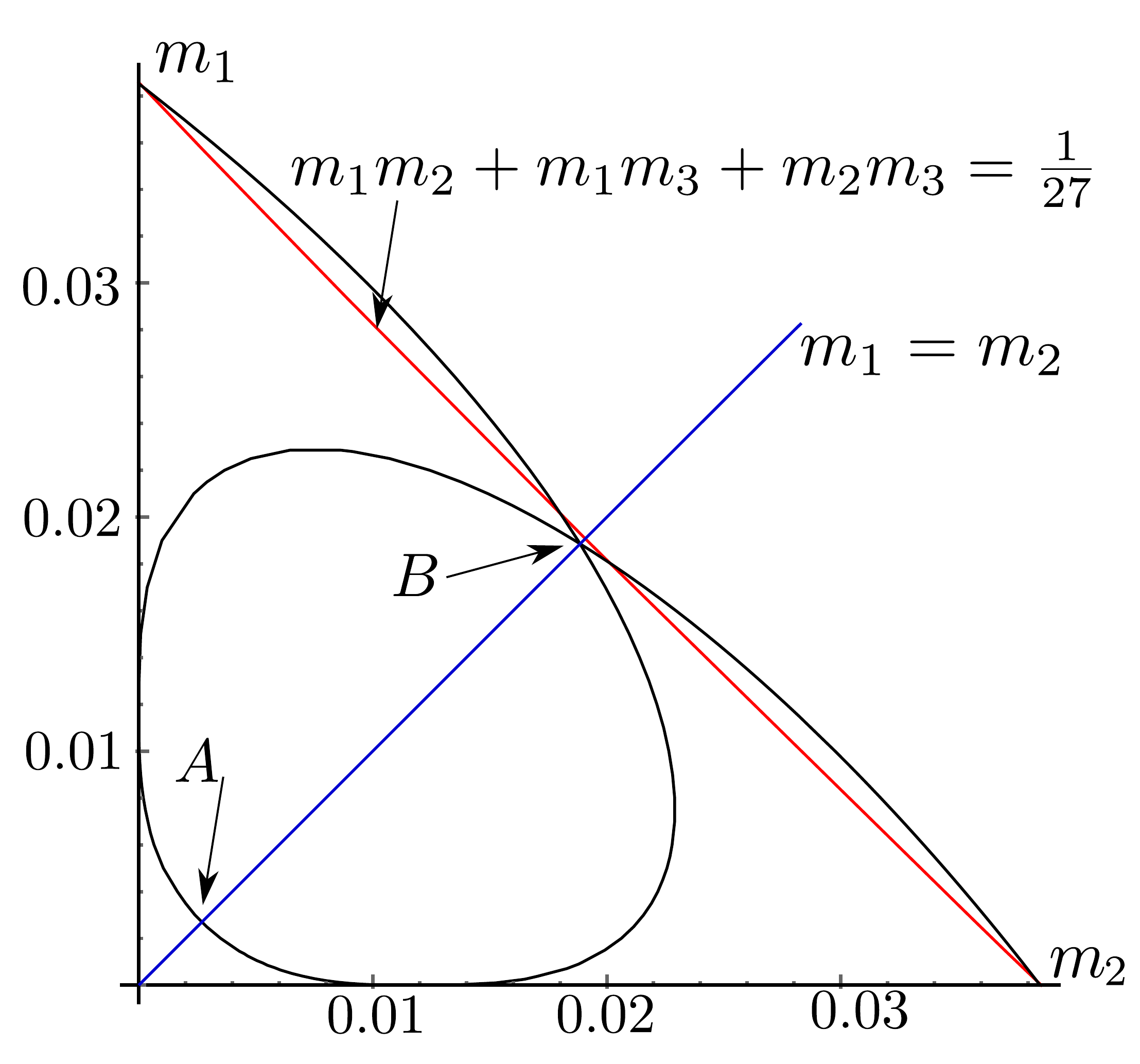}
\caption{Resonance $1$:$1$   with two primary bodies with equal masses in the region II: $A (0.002736,0.002716)$  for $L_5$ and $B(0.01883,0.01883)$ for $L_3$ and $L_6$.}\label{res-1-1sola_inf}
\end{figure}

\begin{figure}[hbt]
\centering
\includegraphics[width=0.5\textwidth]{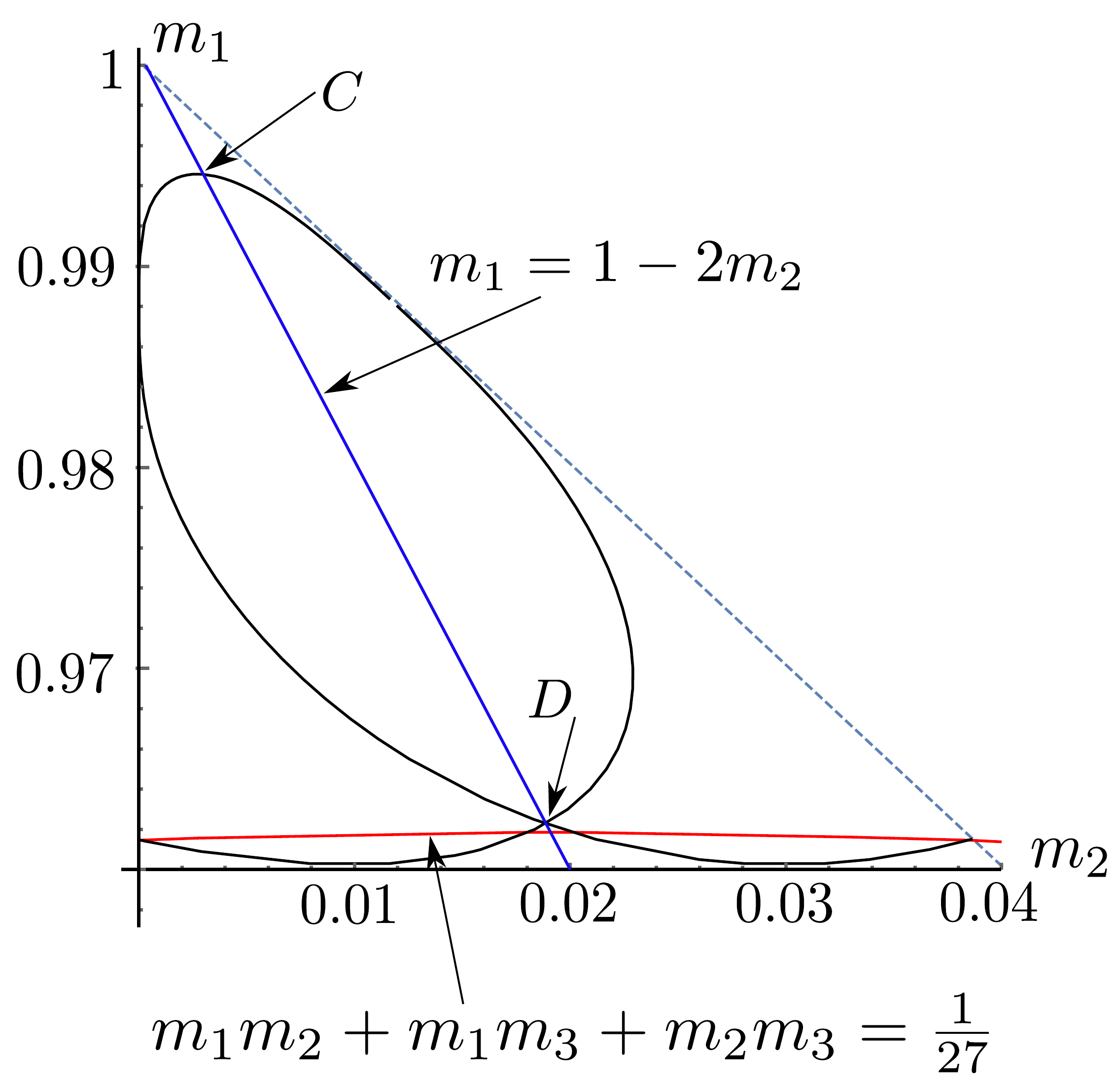} 
\caption{Resonance $1$:$1$   with two primary bodies with equal masses in the region III: $C(0.002736,0.994528)$  for $L_3$ and $D(0.01883,0.96234)$ for $L_5$ and $L_6$.}\label{res-1-1sola_sup}
\end{figure}

\begin{figure}[hbt]
\centering
\subfigure[The different colors indicate the total number of stable equilibrium points in region I. It is from \cite{zotos2020}.]{\includegraphics[width=0.45\textwidth]{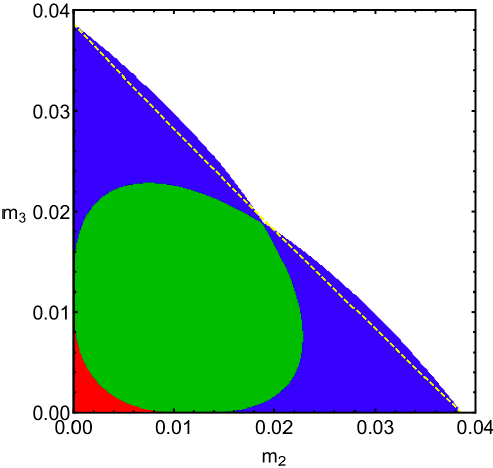}}
\hspace{0.4cm}
\subfigure[Resonance $1$:$1$ curve on region I and  the corresponding stable equilibrium points in the different areas.]{\includegraphics[width=0.5\textwidth]{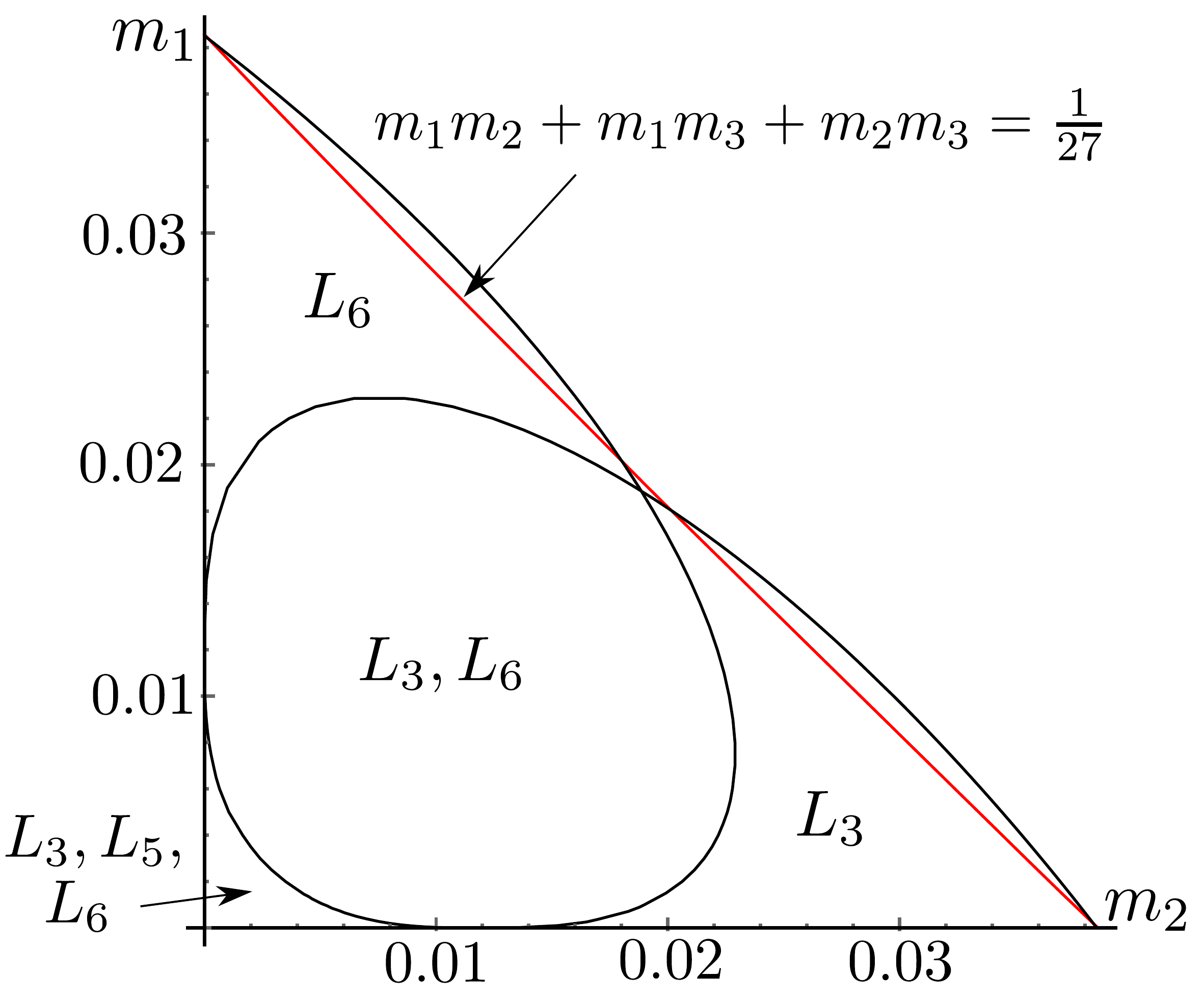}}
\caption{(a) The different colors indicate the total number of stable equilibrium points in region I. (b) Resonance curves $1$:$1$ on region I and stability domain for $L_3$, $L_5$ and $L_6$. For some values of the masses the equilibrium points $L_3$ and $L_6$ are stable even though the Lagrange triangle is not.} \label{fig_cebolla1}
\end{figure}

\begin{figure}[hbt]
\centering
\subfigure[The different colors indicate the total number of stable equilibrium points in region III. It is from \cite{zotos2020}.]{\includegraphics[width=0.47\textwidth]{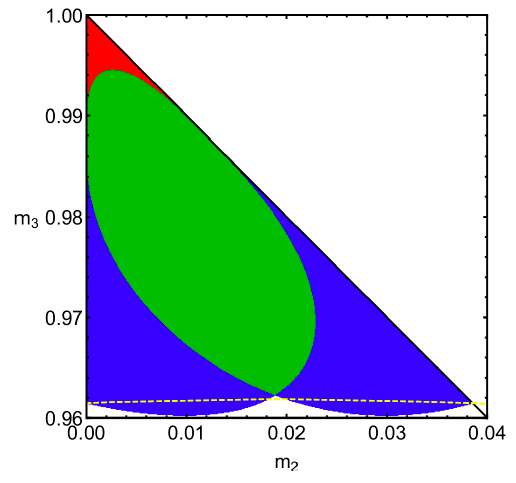}}
\hspace{0.15cm}
\subfigure[Resonance $1$:$1$ curve on region III  and the corresponding stable equilibrium points in the different areas.]{\includegraphics[width=0.46\textwidth]{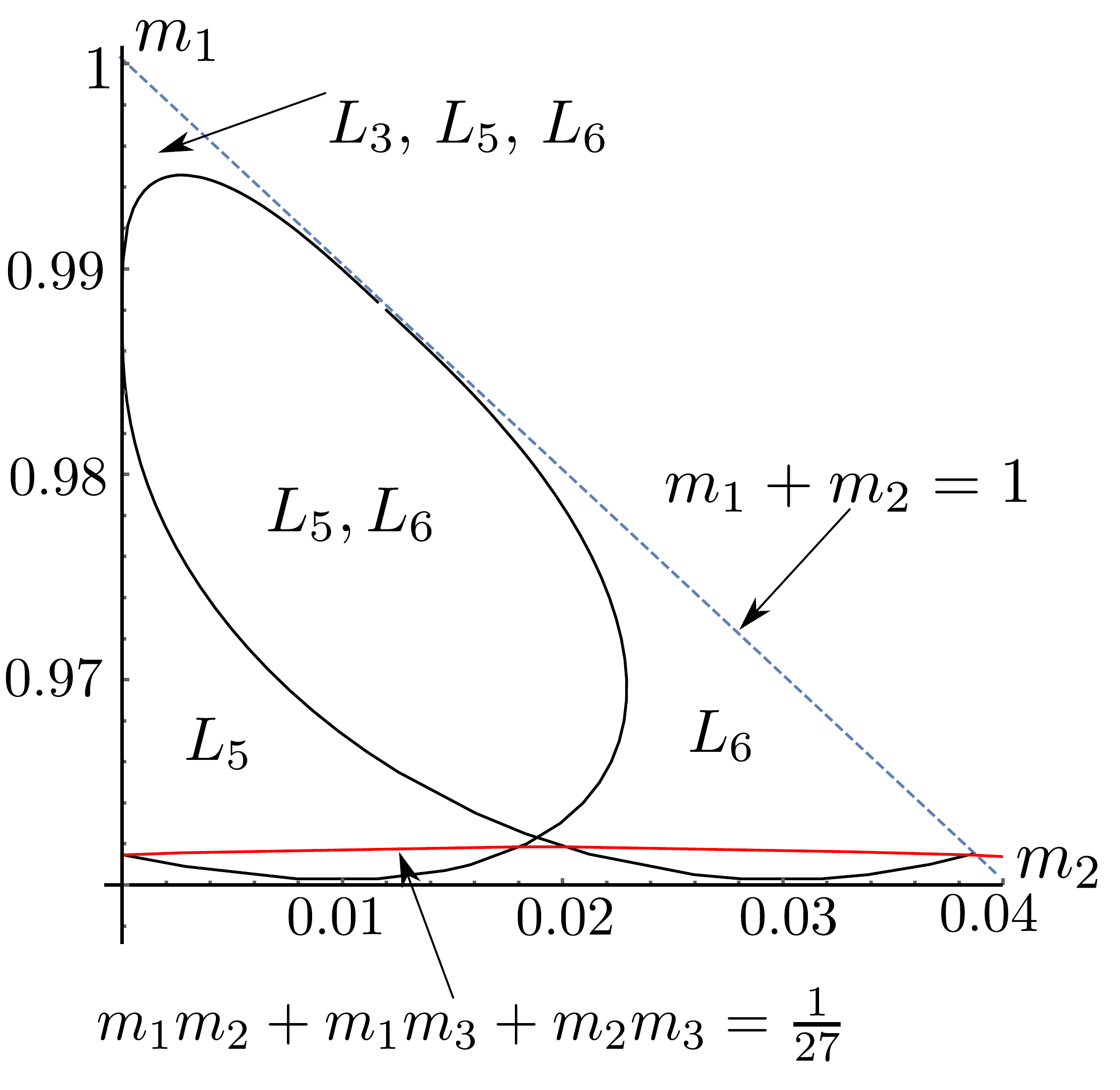}}
\caption{(a) The different colors indicate the total number of stable equilibrium points in region III. 
(b) Resonance $1$:$1$ curve  boundary of the stability region. For some vales of the masses the equilibrium points $L_5$ and $L_6$ are stable even though the Lagrange triangle is not.}\label{fig_cebolla2}
\end{figure}

In Figure \ref{res-1-1sola_inf}, we present  all $1$:$1$ resonance curves in region I. 
On this basis,  we note that there are two points of intersection of
 the resonance curve $1$:$1$ with the straight line $m_1=m_2$ (two primary bodies with equal masses).
The first one is  $A (0.002716,0.002716)$ that  is precisely the point with resonance  $1$:$1$ for $L_5$ where   Burgos and Delgado \cite{burgos-delgado}  ($L_2$ in their notation) established  the existence of a ``blue sky catastrophe", and
the other is $B(0.01883,0.01883)$ resonance values for  $L_3$ and $L_6$,  whose stability has not hitherto been studied.
On the other side, the  intersection point of the resonance curve $2$:$1$ with  line $m_1=m_2$ in Figure \ref{fig4LT}  
corresponds to the equilibrium where  Alvarez-Ram\'{\i}rez et al. \cite{stuchi2015} studied the non-linear stability. 
Actually, similar points are found in region III.  More precisely, we obtain the points by intersecting the resonance curves $1$:$1$ with the line  $m_1=1-2m_2$, that means that $m_2=m_3$, namely  $C(0.002736,0.994528)$  for $L_3$ and $D(0.01883,0.96234)$ for $L_5$ and $L_6$, see Figure \ref{res-1-1sola_sup}. 
A novelty of this article lies in the fact that this is the first time that these points are obtained.
It is interesting to note that in the region II are located the points with $m_1=m_3$ in $1$:$1$ resonance. 

\begin{figure}[hbt]
\subfigure[Domains of stability of relative equilibria $P_{45}$]{\includegraphics[scale=0.3]{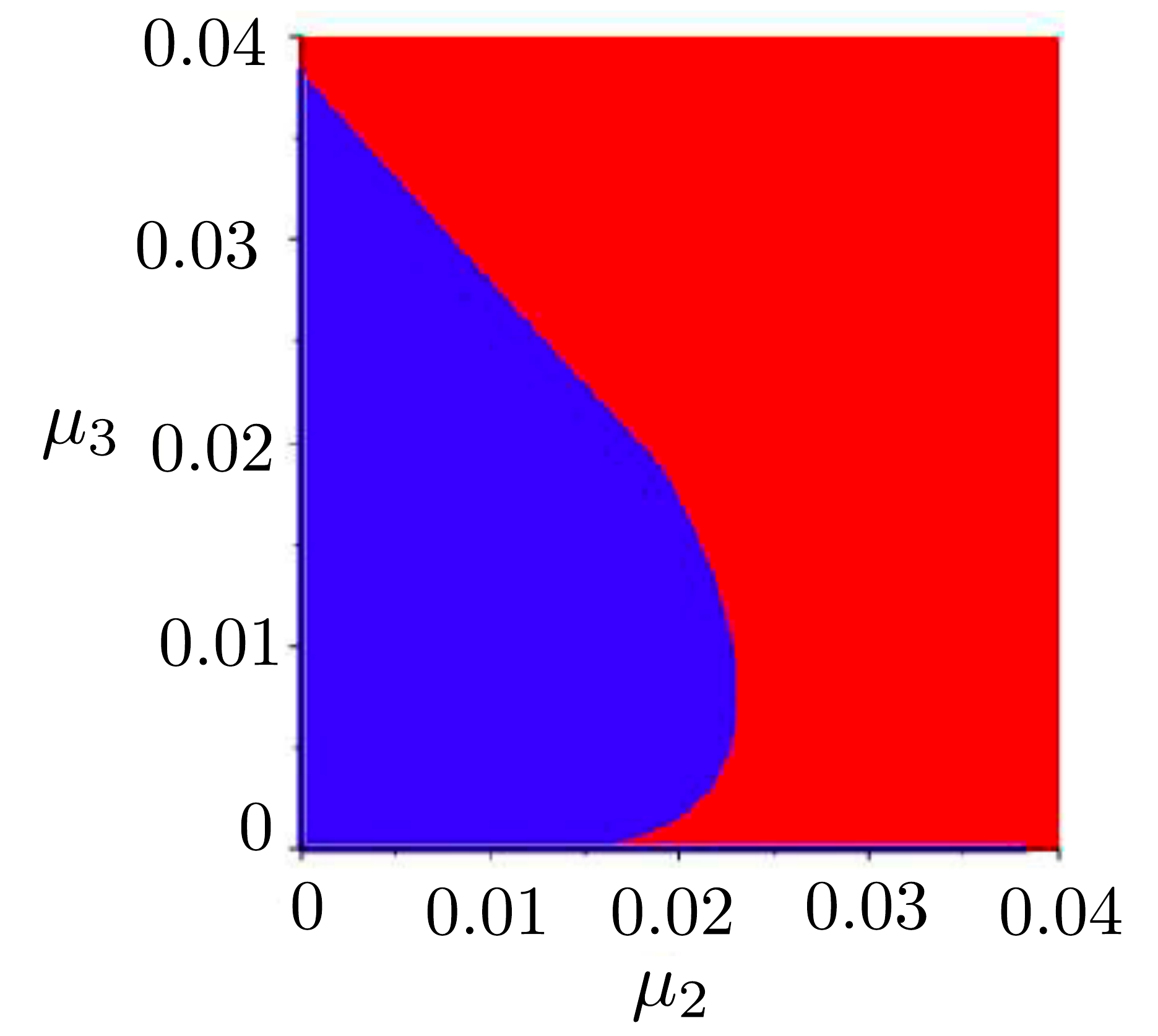}}
	\hspace{0.1cm}
	\subfigure[Domains of stability of relative equilibria $P_{54}$]{\includegraphics[width=0.37\textwidth]{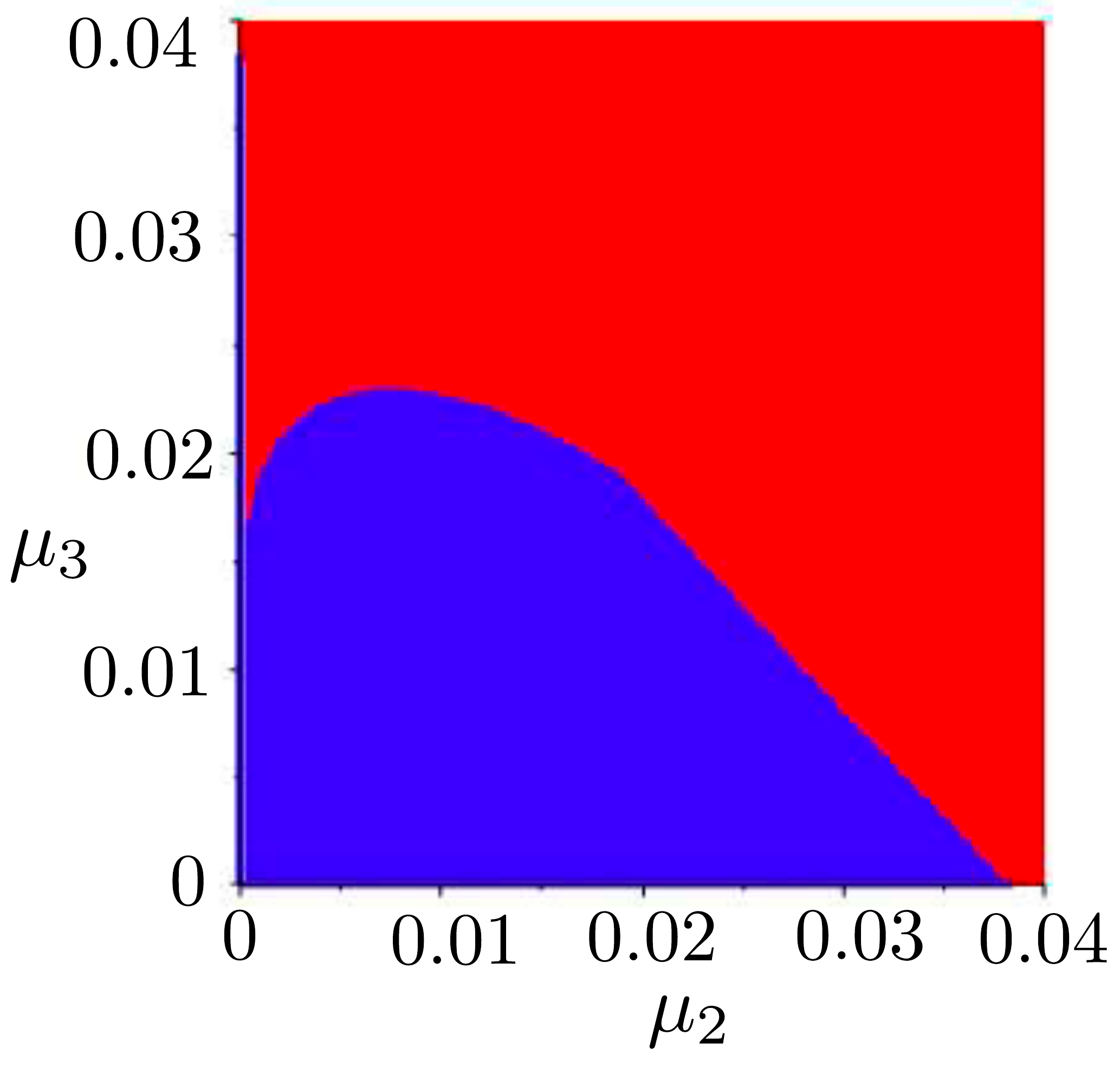}}
\hspace{0.1cm}
	\subfigure[Domains of stability of relative equilibria $P_{55}$]{\includegraphics[width=0.39\textwidth]{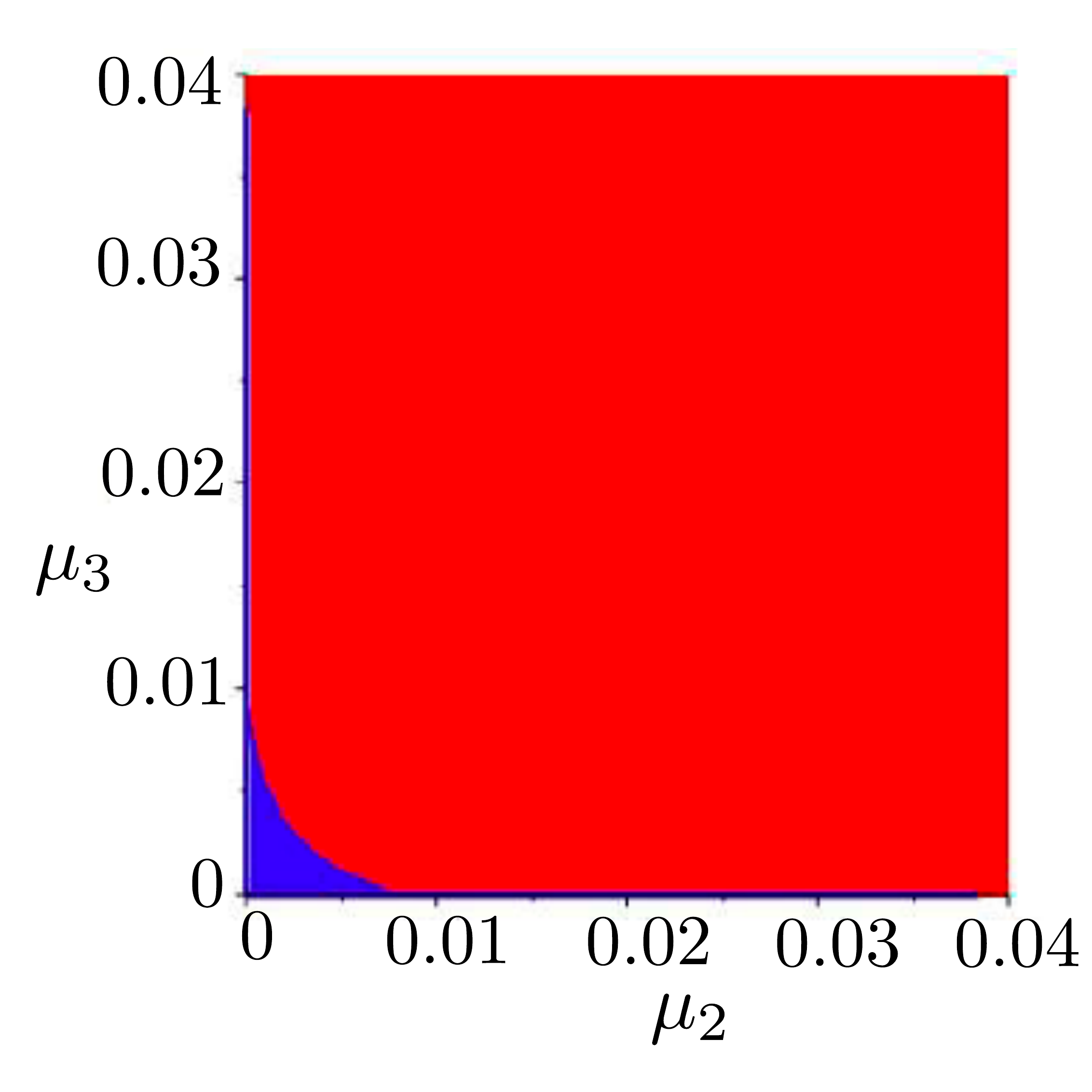}}
	\caption{Graphs computed in \cite{Bardin2021}, where  instability domains are indicated by red color and domains of linear stability  are indicated by blue color.  } \label{fig_bardin}
\end{figure}

Let us quickly mention the results recently reported by Bardin and Volkov \cite{Bardin2021}. They   performed an  analysis of stability and bifurcation of the equilibria 
in the ERFBP  in the case that the Lagrangian triangular configuration of three bodies is stable.
In particular, they  show that the bifurcation is only possible in cases of degeneracy, when the mass of a
primary body vanishes, that is, when the problem degenerates into two copies of the restricted three-body problem.
Bardin and Volkov introduced the mass parameters $\mu_2=\frac{m_2}{m_1+m_2+m_3}$, $\mu_3=\frac{m_3}{m_1+m_2+m_3}$ and denoted by $P_{ij}$  the relative equilibrium that degenerates into a triangular libration point of the restricted three body problem $L_i^{(2)}$ as $\mu_2\to 0$, and degenerates into another triangular  libration point $L_j^{(3)}$ as $\mu_3\to 0$. They claimed that $P_{45}$, $P_{54}$ and $P_{55}$  are the only points that can be both stable and unstable depending on values of parameters. The  stability diagrams were constructed in the plane of parameters $\mu_2$ and $\mu_3$,  and are shown in Figure \ref{fig_bardin}. Hence, comparing the graphs in  this figure with those in Figure   \ref{fig4L3}(a) and Figure \ref{fig4L5}(a),  
we can conclude that their linear stability results are similar to ours, where $P_{45}\to L_6$, $P_{54}\to L_3$
and $P_{55}\to L_5$. It follows that the border between the stable and unstable regions computed by them  are the resonance curves $1$:$1$. 

\begin{figure}[h!]
	\subfigure[It is  Figure 7(a)  taken from \cite{Simo1978}.]{	\includegraphics[width=0.28\textwidth]{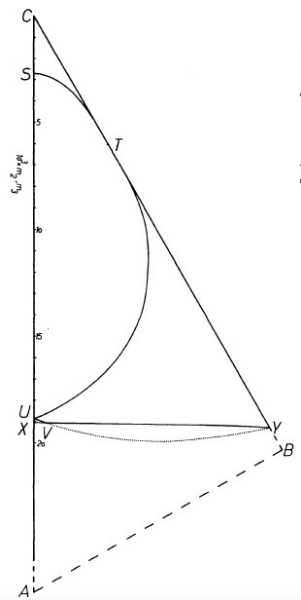}}
	\hspace{1.5cm}
	\subfigure[Figure 7(a) with label on the curves and points obtained by \cite{Simo1978}.]{\includegraphics[width=0.58\textwidth]{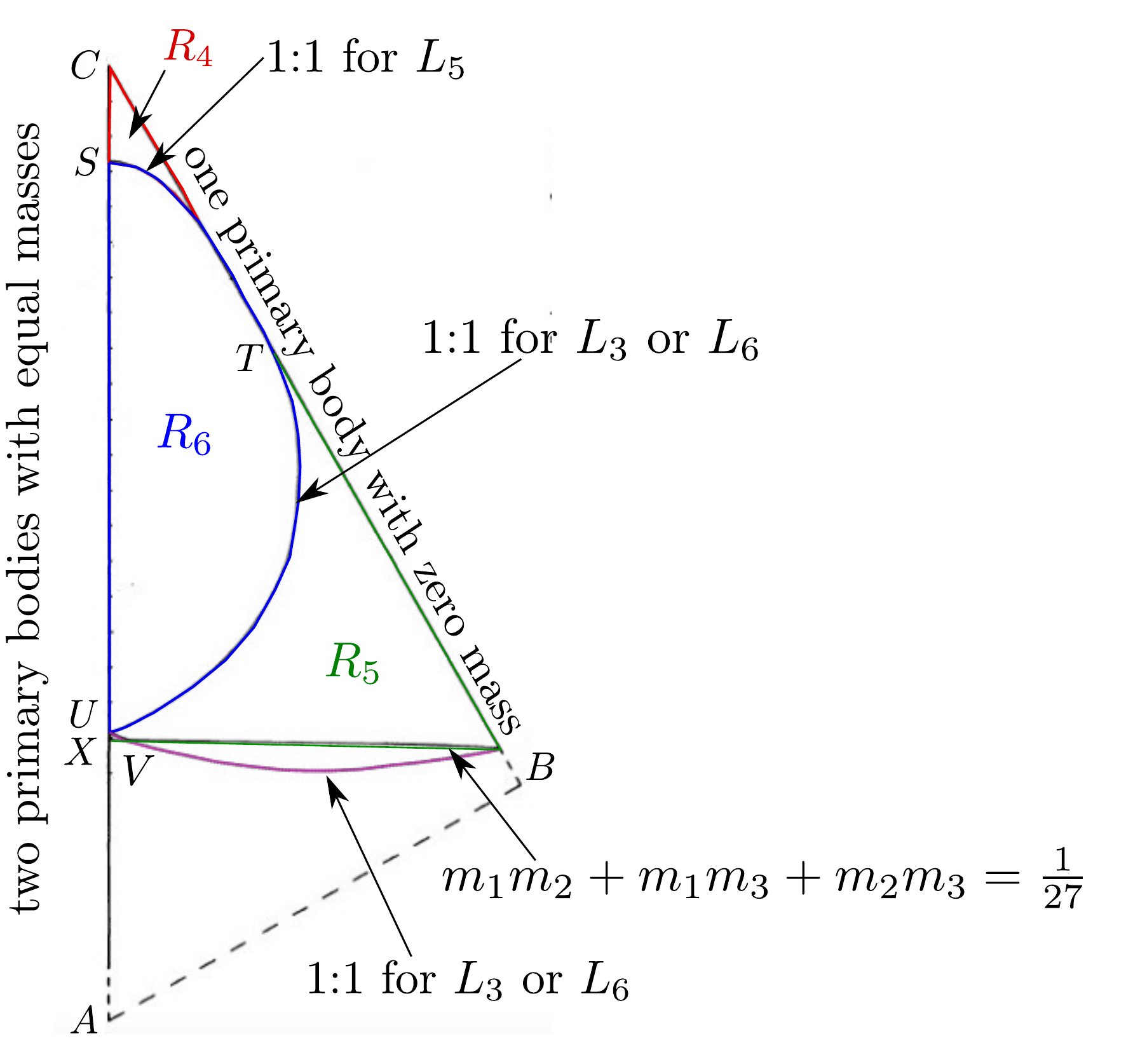}}
\caption{The region $CYVX$ is the linear stability area of the primaries in the mass space, where 
 $R_4=CST$, $R_6=CUT$ and  $R_5=CUVY$ are  stability regions of equilibrium points  $L_5$, $L_3$ and $L_6$.} \label{fig_simo}
\end{figure}

Now, we restrict our attention to the remarkable  paper by   Sim\'o \cite{Simo1978},
where a detailed study is performed by analyzing  numerically the relative equilibrium solutions in the four body  problem.
In particular, we point out in his  Figure 7(a), where   he  displays the  domain  of linear stability of relative equilibria for the ERFBP in $1/6$ of the mass triangle, namely a half of region  I from our Figure \ref{fig_reg1}.  To proceed in the analysis of his results, we note that the coordinates $(m_2,m_3)$ of the points marked in his plot are
 \begin{eqnarray*}
 && A(\frac{1}{3},    \frac{1}{3}), \quad B(0,\frac{1}{2}), \quad C(0,0), \quad S(0.0027096,0.0027096), \quad
\; Y(0,0.038521),   \\
  && T(0,0.011947), \; U(0.018858,0.018858), \;   X(0.019034,0.019064), \; V(0.018114,0.020014).
 \end{eqnarray*}
 His study was confined to the region  where  the three massive bodies configuration is stable, more precisely  at $CYVX$  (see  Figure \ref{fig_simo}). Sim\'o shows  that  the region of stability is separated into three subregions 
 given by  $R_4=CST$, $R_6=CUT$, $R_5=CUVY$, such that $R_4\subset R_6\subset R_5$.
By our calculations,  we are able to identify  the points marked by Sim\'o on the boundary of the region  $CYVX$, and what we get is shown in the Table \ref{tabla_simo} and  displayed in Figure \ref{fig_simo}. We must remember that   $m_1=m_2$ means that 
two primary bodies have equal masses.

The reader should note that the borders of the stability regions determined by Sim\'o are just
the resonance curves $1$:$1$, more precisely,
$ST$  for  $L_5$, $TU$ for $L_6$  and the arc $YVU$ for $L_3$.  Futhermore,  the regions  $R_4\subset R_6\subset R_5$ are 
those shown in the Figures \ref{fig_cebolla1} and \ref{fig_cebolla2}.
For more details see  Figure \ref{fig_simo} where labels are shown.

\begin{table}[hbt] 
\begin{tabular}{|c| c|}
\hline
point on  border of $CYVX$ &  corresponds to  \\ \hline
$U(0.018858,0.018858)$ &   {\scriptsize  point   at  $1$:$1$ resonance curves of $L_3$ and $L_6$ with line $m_1=m_2$:}  \\ 
&  {\scriptsize    point $B$  in   Figure \ref{res-1-1sola_inf} }\\ \hline
$ V(0.018114,0.020014)$  &  {\scriptsize intersection point between $1$:$1$ resonance curve of $L_3$} \\
 &  {\scriptsize and Routh's critical curve:  point $B_{L_3}$ in Figure \ref{fig4L3} (a)} \\ \hline 
 $T(0,0.011947)$ &  {\scriptsize intersection point between $1$:$1$ resonance curve of $L_5$} \\
  & {\scriptsize either  $m_1=0$ or $m_2=0$: Figure \ref{fig4L5} (a)}\\ \hline
  $S(0.0027096,0.0027096)$ & {\scriptsize  point   at  $1$:$1$ resonance curve of $L_5$ with line $m_1=m_2$:}  \\ 
&  {\scriptsize   point $A$  in   Figure \ref{res-1-1sola_inf}}\\ \hline
$X(0.019034,0.019064)$ & {\scriptsize intersection point between the line $m_1=m_2$}\\
 &  {\scriptsize and Routh's critical curve: Figure  \ref{fig4L6}}\\ \hline
 $Y(0,0.038521)$ &  {\scriptsize intersection point between the Routh's critical curve and }\\
 & {\scriptsize $1$:$1$ resonance curve of $L_6$ either $m_1=0$ or $m_2=0$:}\\ 
 &   {\scriptsize points $R_{L_6}$ and ${\mathcal R}_{L_6}$  in Figure  \ref{fig4L6} } \\\hline
\end{tabular}\label{tabla_simo}
\caption{Relationship between the values obtained by Sim\'o those obtained by us. }
\end{table}

\section{Concluding remarks} \label{sec5}
This paper summarizes the most  known results (up to this date) about the location, counting and linear stability of the equilibrium points in the ERFBP.

In this review, the main attention was paid in the unification of known results about relative equilibria or  central configuration of the planar restricted $(3 + 1)$-body problem  with primaries in Lagrange's configuration.

We would like to highlight  that resonance curves are of utmost importance  in determining the linear stability domain for the equilibrium points.
Some resonance curves in region I  were previously calculated  by Budzko and Prokopenya, see \cite{Budzko2009} and \cite{BudzkoProkopenya}. 
However, we have given a step forward  by computing some other resonance curves by taking into consideration the three regions in the plane where the Lagrange's configuration turn out to be stable. This allowed us to  determine regions, in a numerical way, in the mass space where the points of equilibrium  are linearly stable, and to find out the regions where these points of equilibrium corresponds.   
Based on  numerical methods, Zotos \cite{zotos2020}   showed the exact number of stable equilibrium points inside the  regions I, II and III. In fact, he claims that  regions where linearly stable points of equilibrium exist almost coincide with the regions for which the Lagrange's configuration of the primaries is stable.
We recall that similar results were reported by Budzco \cite{Budzko2009}, but  only for  the equilibrium point $S_1$ ($L_3$ in our notation, see Table \ref{tabb}).
Our analysis provide firm numerical evidence that  $L_3$ and $L_6$ in region I, as well as $L_5$ and $L_6$ in regions II and III are  stable in linear approximation for some mass values $m_1$, $m_2$ and $m_3$ which are outside the stability domain of the Lagrange's triangle. Furthermore, 
we present numerical results to show that the  boundaries of the stability domain are determined  by the $1$:$1$ resonance curves.
  At present, this is a remarkable fact for which we have no explanation yet to offer.  This is left, for now, as a future avenue of research on the ERFBP.

\section*{Acknowledgements}
The first author was supported by a  UAM fellowship of doctoral studies.  The second author  was
partially supported by Special Program to Support Teaching and Research Projects 2021 from  CBI UAM-Iztapalapa   (M\'exico).

The authors would like to thank the anonymous reviewers for their comments and fruitful suggestions to improve the quality of the paper.

\section*{Statements and Declarations}  

The authors declare that they have no conflict of interest. This study was funded by the Universidad Aut\'onoma Metropolitana  (Metropolitan Autonomous University) Mexico.

Jos\'e Alejandro Zepeda Ram\'{\i}rez and Martha Alvarez-Ram\'{\i}rez   contributed to the design and implementation of the research, to the analysis of the results and to the writing of the manuscript.

Data sharing is not applicable to this article as no new data were created or analyzed in this study.

\end{document}